\documentclass[12pt]{article}

\usepackage{url}
\usepackage{fancyhdr}
\usepackage{setspace} 
\usepackage{pdfscreen} 
\usepackage{mathptmx}
\usepackage{latexsym}
\usepackage{amssymb}
\usepackage{amsthm}
\usepackage{amsmath}
\usepackage{authblk}

\usepackage{graphicx}

\begin{document}

\date{}
\title{Determination of the 4-genus \\ of a complete graph \\ (with an appendix)}

\author{\Large{Serge Lawrencenko, Beifang Chen, Hui Yang}}

\affil[]{}
\renewcommand\Authands{}
\maketitle

\noindent {\bf Abstract.} In this paper, the quadrangular genus (4-genus) of the complete graph  $K_p$  
is shown to be  $\gamma_4 (K_p) = \lceil {p(p-5)}/{8} \rceil +1$
for orientable surfaces. This means that $K_p$ is minimally embeddable in the closed orientable surface of genus $\gamma_4 (K_p)$ under the constraint that each face has length at least 4. In the most general setting, the genus of the complete graph was established by Ringel and Youngs and was mainly concerned with triangulations of surfaces. Nonetheless, since then a great deal of interest has also been generated in quadrangulations of surfaces. Hartsfield and Ringel were the first who considered minimal quadrangulations of surfaces. 

Sections 1--4 of this paper are essentially a reproduction of the original 1998 version as follows:
Chen~B., Lawrencenko~S., Yang~H. Determination of the 4-genus of a complete graph, submitted to Discrete Mathematics and withdrawn by S. Lawrencenko, June 1998, URL: {\url{https://t.co/cUg6R9Jwyw}}.
More discussion on this 1998 version is held and some copyright issues around the quadrangular genus of complete graphs are clarified in the Appendix to the current version of the paper; the Appendix was written in 2017.

\medskip

\noindent {\bf Keywords:} graph embedding, genus of graph
\medskip

\section {Introduction}

In the most general setting, the genus of the complete graph was established by Ringel and Youngs [3] and was mainly concerned with {\it triangulations} of surfaces. Nonetheless, since then a great deal of interest has also been generated in {\it quadrangulations} of surfaces. In particular, Hartsfield and Ringel [2] were the first who considered minimal quadrangulations of surfaces; see the Appendix for further information and references. 
Sections 1--4 of this paper are essentially a reproduction of the original 1998 version as follows:
Chen~B., Lawrencenko~S., Yang~H. Determination of the 4-genus of a complete graph, submitted to Discrete Mathematics and withdrawn by S. Lawrencenko, June 1998. URL: {\url{https://t.co/cUg6R9Jwyw}}.
More discussion on this 1998 version is held and some copyright issues around the quadrangular genus of complete graphs are clarified in the Appendix to the current version of the paper; the Appendix was written in 2017~\cite{L4}.

Given an embedding of a graph in a surface, the number of edges on the boundary of a face is called the {\it length of the face.} We define the 
{\it $n$-genus of a graph} $G$, denoted $\gamma_n (G)$, to be the minimum genus $h$ of orientable surface $S_h$ in which $G$ is 2-cell embeddable with each face having length at least $n$. Given integers $p$ and $n$ such that $3 \leqslant n \leqslant p$, the general problem is to evaluate the $n$-genus of the complete graph $K_p$. The celebrated formula for the genus of the complete graph obtained by Ringel and Youngs [3] answers this question for $n = 3$, $p \geqslant 3$. Furthermore, Hartsfield and Ringel [2] constructed quadrangular embeddings of $K_p$ for $p \equiv 5 \mod 8$, which partly answers the question for $n = 4$ and, in this paper, we will focus on this case. Our main result is the 4-genus (or the {\it quadrangular genus}) of the complete graph:

\medskip
\noindent
\hbox{\bfseries Theorem 1.}
{\itshape
$$
\gamma_4 (K_p) = \bigg{\lceil} \frac {p(p-5)}{8} \bigg{\rceil} +1.
$$
}
\medskip

We review some related results on the subject in Section 2, construct quadrangular embeddings of specific graphs in surfaces in Section 3, and prove Theorem 1 in Section 4. The terminology and notation are standard and will be found to conform, for the most part, for instance, with [4, 5].

\section {Current graphs}

We assume the familiarity with the rotation scheme technique; we refer the reader to [4, 5] for details. In this section, for the sake of completeness, we review the main points.

Let $G$ be a given graph. To construct 2-cell embeddings of $G$  in orientable surfaces, the following technique is used. Let $\Gamma $  be a group of order $n$  whose neutral element will be denoted by $0$ and let $\Omega $ be a subgroup of $\Gamma$. Let $K$ be a pseudograph and let 
$K^*=\{(u,v): uv \in E(K)\}$ denote the set of its oriented edges. A {\it current graph} for $\Gamma$ is a quadriple $(K, \Gamma, \varphi, \lambda )$. Here $\varphi$ is a rotation scheme for $K$, $\lambda$ is an assignment of ``currents'' (which are associated with the elements of $\Gamma$) to the oriented edges of $K$, which is defined to be a mapping $\lambda$: 
$K^* \rightarrow \Gamma \setminus \{ 0 \}$ such that $\lambda (v, u) = - \lambda (u, v)$  for each 
$(u, v) \in K^*$.  The index of the subgroup $\Omega$ in $\Gamma$ is called the {\it index of the embedding}, or the {\it index of the current graph}, denoted $[ \Gamma : \Omega ]$. Each right coset of $\Omega $ in $\Gamma $  determines one row of the corresponding rotation scheme for $G$,  and these $[\Gamma : \Omega]$  rows generate the whole rotation scheme, which corresponds to a 2-cell embedding of $G$  in some orientable surface.

Denote by $\pi_v$ the {\it rotation} at vertex $v$,  that is, the permutation of the vertices adjacent to $v$, given by the rotation scheme. In what follows, we shall assume that addition is always carried out in the group $\Gamma$. We say that the Kirchhoff current law (briefly, KCL) holds at vertex $v \in V(K)$  if the sum $\sum_v$ of the currents directed away from $v$ is $0$. 

A 2-cell embedding of $G$ is said to be {\it $n$-angular} if all its faces have length $n$.  The following is a collection of sufficient conditions which in combination guarantee that the current graph determines an $n$-angular embedding of the complete graph $K_p$ in some orientable surface. Here we address the index 1 current graphs. This is the simplest case, in which $\Gamma = \Omega = {\mathbb Z}_p$  is the additive group of integers modulo $p$.

\medskip
\noindent
\hbox{\bfseries Lemma 1.}
{\itshape
An index 1 current graph $(K, {\mathbb Z}_p, \varphi, \lambda)$ determines an $n$-angular embedding of the complete graph $K_p$ in an orientable surface if it satisfies the following conditions:
\medskip

\noindent {\rm {(C1)}} $\varphi$ induces a 2-cell embedding of $K$ in $S_h$ having exactly one 2-cell; 

\noindent {\rm {(C2)}} this 2-cell is a $(p - 1)$-gon;

\noindent {\rm {(C3)}} each element of ${\mathbb Z}_p \setminus \{0\}$ appears exactly once as a current on some oriented edge of the $(p - 1)$-gon (traced clockwise, for certainty);

\noindent {\rm {(C4)}} each vertex of $K$ has valence $n$ ($n \geqslant 3$);

\noindent {\rm {(C5)}} the KCL holds at each vertex.
}
\medskip
\begin{proof}
Label the verticies of $K_p$ with the elements of the group ${\mathbb Z}_p = \{ 0, 1, 2, \ldots, {p-1} \}$. By conditions C1, C2, and C3, the clockwise-oriented boundary of the 2-cell can serve as row $0$ for the rotation scheme:
 
$$
\pi_0: (a_1, a_2, \ldots , a_{p-1}).
$$

\noindent We then obtain row $i$, $1 \leqslant i \leqslant p-1$, by adding $i$ to each entry in row $0$:

$$
\pi_i : (a_1 + i, a_2 + i, \ldots, a_{p-1}+ i).
$$

    \begin{center}
    \hspace*{0cm}    \pdfimage width 0.5\textwidth {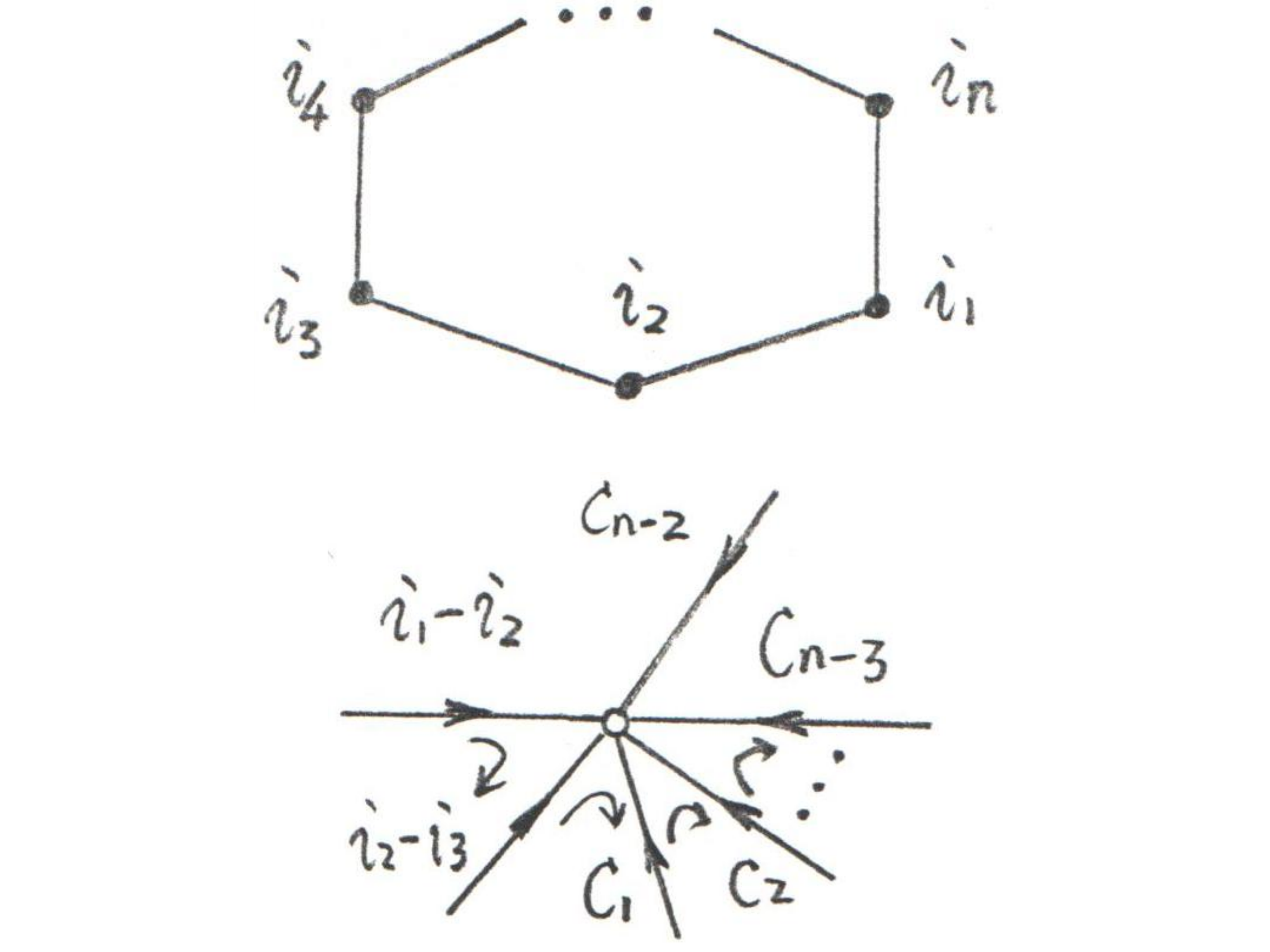}
\newline

    Figure 1.
    \end{center}

\noindent We next use conditions C4 and C5 to show that the 2-cell embedding determined by this rotation scheme is $n$-angular. For this, consider two successive boundary edges $(i_2, i_1)$ and 
$(i_2, i_3)$ for an arbitrary face, as depicted in Figure 1 (top). Necessarily, $\pi_{i_2} (i_1) = i_3$,  so that 
$\pi_0 (i_1 - i_2) = i_3 - i_2$.  Thus, in the current graph $K$ we have the situation depicted in 
Figure~1 (bottom). Then, by C5, we have:

$$
a_1 + a_2 + a_3 + a_4 + \cdots + a_{n-2} + i_1 - i_3 = 0.
$$

\noindent Now, let 
$$
i_4 = i_3 - c_1, \quad i_5 = i_4 - c_2, \quad \ldots,  \quad i_{n-1} = i_{n-2} - c_{n-4},
$$
$$
i_n = i_{n-1} - c_{n-3}, \quad  i_{n + 1} = i_n - c_{n-2}.
$$

\smallskip

\noindent Finally, condition C5 implies that $i_{n + 1} = i_1$.  Thus, all faces have length $n$.
\end{proof}

Especially for $n = 4$,  it is useful to observe that verification of whether conditions C1--C5 hold is equivalent to checking the following ``quadrangular rule'': 

If in row $i$ we find ``$\ldots, \, j, \, k, \, \ldots$''

and in row $k$ we find ``$\ldots, \, i, \, l, \, \ldots$'', 

then in row $j$ we find ``$\ldots, \, l, \, i, \, \ldots $''

and in row $l$ we find ``$\ldots, \, k, \, j, \, \ldots $''.
    \begin{center}
    \hspace*{0cm}    \pdfimage width 0.4\textwidth {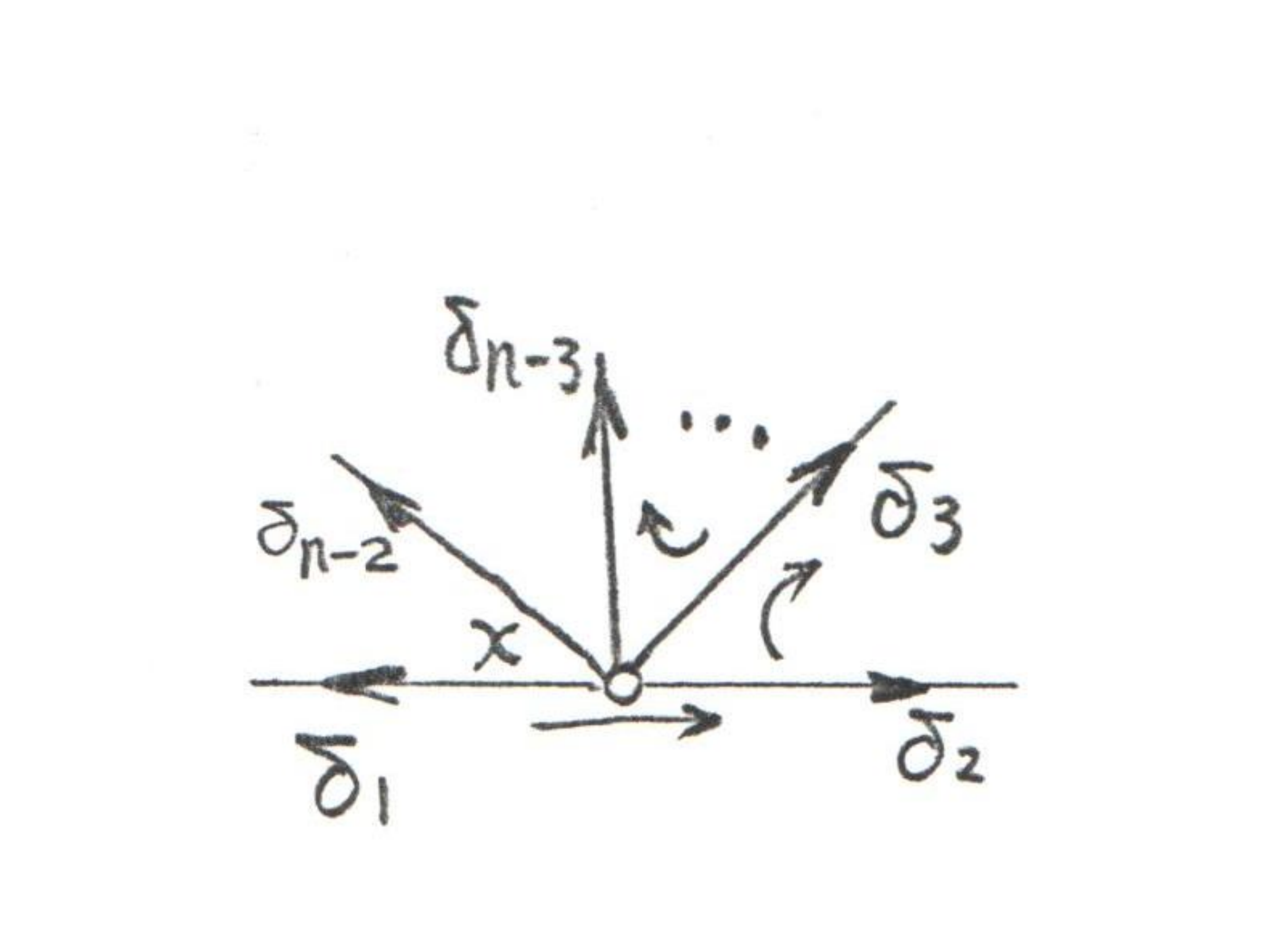}
 \end{center}
    \begin{center}
Figure 2.
    \end{center}

\bigskip

Lemma 1 generalizes the method employed by Ringel and Youngs [3] for construction of triangular embeddings of $K_p$. To obtain $n$-angular embeddings, sometimes higher index graphs may be required. In [3], current graphs (after some appropriate modifications) are also used to construct triangular embeddings of nearly complete graphs $K_m + {\overline {K}_{p - m}}$, where ``$+$'' denotes the join of the summands and $m$ is near $p$. Similarly, in order to determine the $n$-genus of $K_p$, we will also need to construct $n$-angular embeddings of nearly complete graphs by using modified current graphs. More precisely, in the index $1$ current graph $(K, {\mathbb Z}_m, \varphi, \lambda)$, condition C4 is modified to allow vertices of valence $n - 2$ so that the sum of the incident currents directed away from any such vertex is of order $m$ as an element of ${\mathbb Z}_m$.  In Figure 2, such vertices are labeled by letters $x$, $y$, $z$, $w$,  etc., and $\sum_i \delta_i$  is a generator of ${\mathbb Z}_m$. More precisely, we consider an index $1$ current graph $(K, {\mathbb Z}_m, \lambda, \varphi)$ satisfying the following conditions:

\medskip
\noindent (C1$^ {\prime}$) The vertex set of $K$ consists of two parts: unlabeled vertices and labeled vertices $x$, $y$, etc.;

\noindent (C2$^{\prime}$) $\varphi$ induces an embedding of $K$ in some orientable surface $S_h$, having exactly one $2$-cell, and this $2$-cell is a $(p - 1)$-gon;

\noindent (C3$^{\prime}$) each element of ${\mathbb Z}_m \setminus \{ 0 \}$ appears exactly once as a current on some oriented edge of the $(m - 1)$-gon (traced clockwise, for certainty);

\noindent (C4$^{\prime}$) each unlabeled vertex has valence $n$ and each labeled one has valence $n - 2$;

\noindent (C5$^{\prime}$) the KCL holds at each unlabeled vertex;

\noindent (C6$^{\prime}$) for each labeled vertex $x$,  the sum of the incident currents directed away from $x$  is an element of order $m$ in ${\mathbb Z}_m$. 

\medskip
Now we construct a rotation scheme for a current graph satisfying C1$^{\prime}$--C6$^{\prime}$. Similarly to the proof of Lemma 1, from the current graph we derive a 2-cell of this form:
$$
(\ldots, \, \, a, \, \, - \delta_1, \, \, \delta_2,  \, \, b, \, \, \ldots, \, \, - \delta_2, \, \, \delta_3, \, \, \ldots, \, \,  -\delta_{n-2}, \, \, \delta_1, \, \, \ldots).
$$

\noindent Then for $i \in {\mathbb Z}_m$ we find

$$
\pi_i : \, \, (\ldots, a + i, \, \, -\delta_1 + i, \, \, x, \, \, \delta_2 + i, \, \, b + i, \ldots, \, \, -\delta_2 + i, \, \,  \delta_3 + i, \ldots, \, \, -\delta_{n-2} + i, \, \, \delta_1 + i, \, \ldots).
$$

\noindent In this fashion we generate $\{ \pi_i  \} _{ i \in \mathbb Z_m}$.  For each labeled vertex $x$, we append

$$
\pi_x : \, \, ((m - 1) \sum \delta_i, \, \, (m - 2) \sum \delta_i, \, \, \ldots, \, \, 2 \sum \delta_i, \, \, \sum \delta_i, \, \, 0)
$$

\noindent to the rotation scheme. If we have $p - m$ labeled vertices $\{ x_i \}$ ($1 \leqslant i \leqslant p- m$) of this type, then the collection

$$
\{ \pi _i  \} _{ i \in \mathbb Z_m} \cup \{ \pi_{x_i} \}_ {i = 1} ^{p-m}
$$

\noindent determines a 2-cell embedding of the join $K_m + \overline K_{p-m}$ in some orientable surface.

\medskip
\noindent
\hbox{\bfseries Lemma 2.}
{\itshape
The embedding constructed above is an $n$-angular embedding.
}
\medskip
\begin{proof}
By Lemma 1, we only need to consider all pairs of neighboring elements in 
$-\delta_1, \, \, x, \, \, \delta_2+i, \, \, -\delta_{n-2}, \, \, \delta_1$ as well as the pairs 
$-\delta_{ j-1}, \, \, \delta_j$  ($j \, = \, 3, \, 4, \, ... , \, n - 2$)  in row 
$\pi_i$  and all successive pairs of elements in row $\pi_x$.  Consider two successive boundary edges 
$(i, \, x)$  and 
$(i, \, \delta_2 + i)$.  Thus, $\pi_i (x) = \delta_2 + i$.  By construction of the rotation scheme, we have:
\begin{equation*}
\pi_{\delta_2 +i}(i) = \delta_3 + \delta_2 + i,
\end{equation*}
\begin{equation*}
\pi_{\delta_3+\delta_2+i}(\delta_2+i) = \delta_4 + \delta_3 + \delta_2 + i,
\end{equation*}
\begin{equation*}
\ldots,
\end{equation*}
\begin{equation*}
\pi_{\delta_{n-2} + \delta_{n-3} + \cdots +\delta_1 + i}
(\delta_{n-2}+\delta_{n-3}+\cdots+\delta_2+i)=x,
\end{equation*}

\noindent and

\begin{equation*}
\pi_x (\delta_{n-2} + \delta_{n-3} + \cdots + \delta_1+i) = i.
\end{equation*}

\noindent So this face has length $n$. Check with Figure 3. 
\begin{center}
\hspace*{0cm}    \pdfimage width 0.6\textwidth {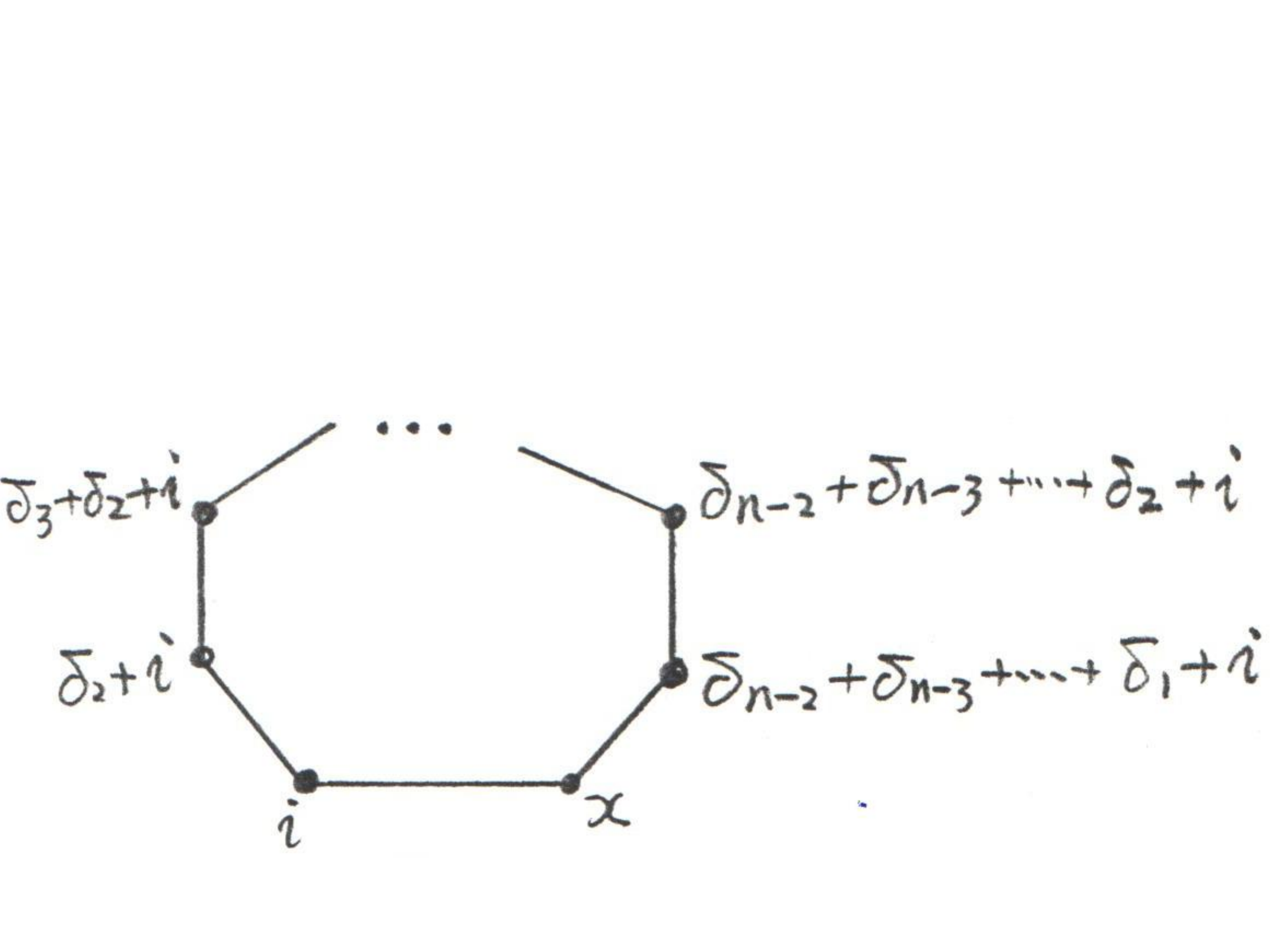}
\end{center}
\begin{center}
    Figure 3.
\end{center}

\bigskip

Similarly, it can be checked that all other faces have length $n$.
\end{proof}
\medskip

(For $n = 3$, the labeled vertices are the ``vortices'' introduced by Ringel and Youngs to construct triangular embeddings of nearly complete graphs.)

Lemmas 1 and 2 suggest a general method to construct $n$-angular embeddings. However, the problem remains to find suitable current assignments.

\section{Quadrangular embeddings}

In this section, we apply current graph theory to construct quadrangular embeddings of complete and nearly complete graphs. By $K_p-K_i$ we denote the graph obtained from the complete graph $K_p$ by deleting all edges of a fixed complete subgraph $K_i$. We will depict the rotation scheme $\varphi$ using black dots $\bullet$ and hollow dots $\circ$; namely, clockwise at $\bullet$ and counterclockwise at $\circ$.

\medskip
\noindent
\hbox{\bfseries Lemma 3.}
{\itshape
There exist quadrangular embeddings of $K_{8s}$  in  $S_{8s^2-5s+1}$ \, for  $s \geqslant 1$.
}
\smallskip

\begin{center}
\hspace*{0cm}    \pdfimage width 1.0\textwidth {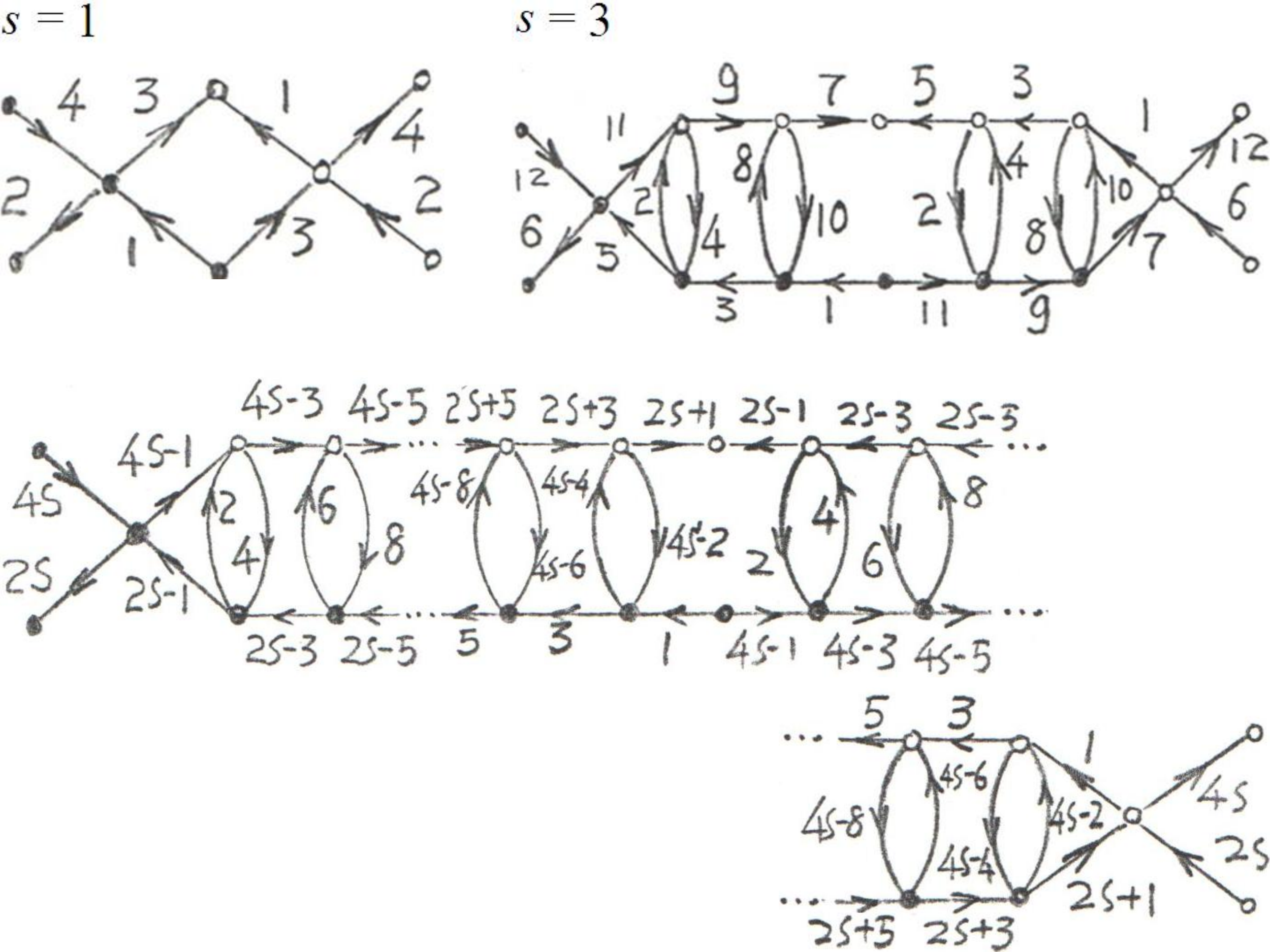}
\end{center}
\begin{center}
    Figure 4.
\end{center}

\begin{proof}
In order to construct quadrangular embeddings of $K_{8s}$, we need index 2 current graphs $(K, {\mathbb {Z}}_{8s}, \varphi, \lambda)$ satisfying the following conditions:
\medskip

\noindent (F1) $\varphi$ induces a 2-cell embedding of $K$ in $S_h$, having exactly two 2-cells;

\noindent (F2) the two 2-cells are both $(8s - 1)$-gons;

\noindent (F3) each element of ${\mathbb Z}_{8s} \setminus \{ 0 \}$ appears exactly once as a current on some oriented edge of each of the $(8s - 1)$-gons (traced clockwise, for certainty);

\noindent (F4) each vertex of $K$ has valence 4, 2, or 1;

\noindent (F5) KCL holds at each vertex of valence 4;

\noindent (F6) for a vertex of valence 1, the incident current has valence 4 or 2 in ${\mathbb Z}_{8s}$, that is, the incident current is $2s$ or $4s$ (respectively), and such currents occur in this order:

$$
 \ldots, 2s, -2s, 4s, \ldots 
$$
or
$$ 
\ldots, 4s, -2s, 2s, \ldots;
$$
(F7) for a vertex $v$ of valence 2, the incident currents belong to
$$
A = \{1, 3, 5, \ldots, 8s - 1 \}
$$
\noindent and the sum $\sum_v$ of incident currents directed away from $v$ is of order 2.
\medskip

The entire rotation scheme is obtained in the following way: Firstly obtain the two rows, $\langle 0 \rangle$ and $\langle 1 \rangle$, from the two clockwise 2-cells. Then obtain row $i$ by adding $i$ to every entry of row $\langle 0 \rangle$ for $i \equiv 0 \mod 2$, and to every entry of row $\langle 1 \rangle$ for  $i \equiv 1 \mod 2$.

\bigskip
\begin{center}
\hspace*{0cm}    \pdfimage width 1.0\textwidth {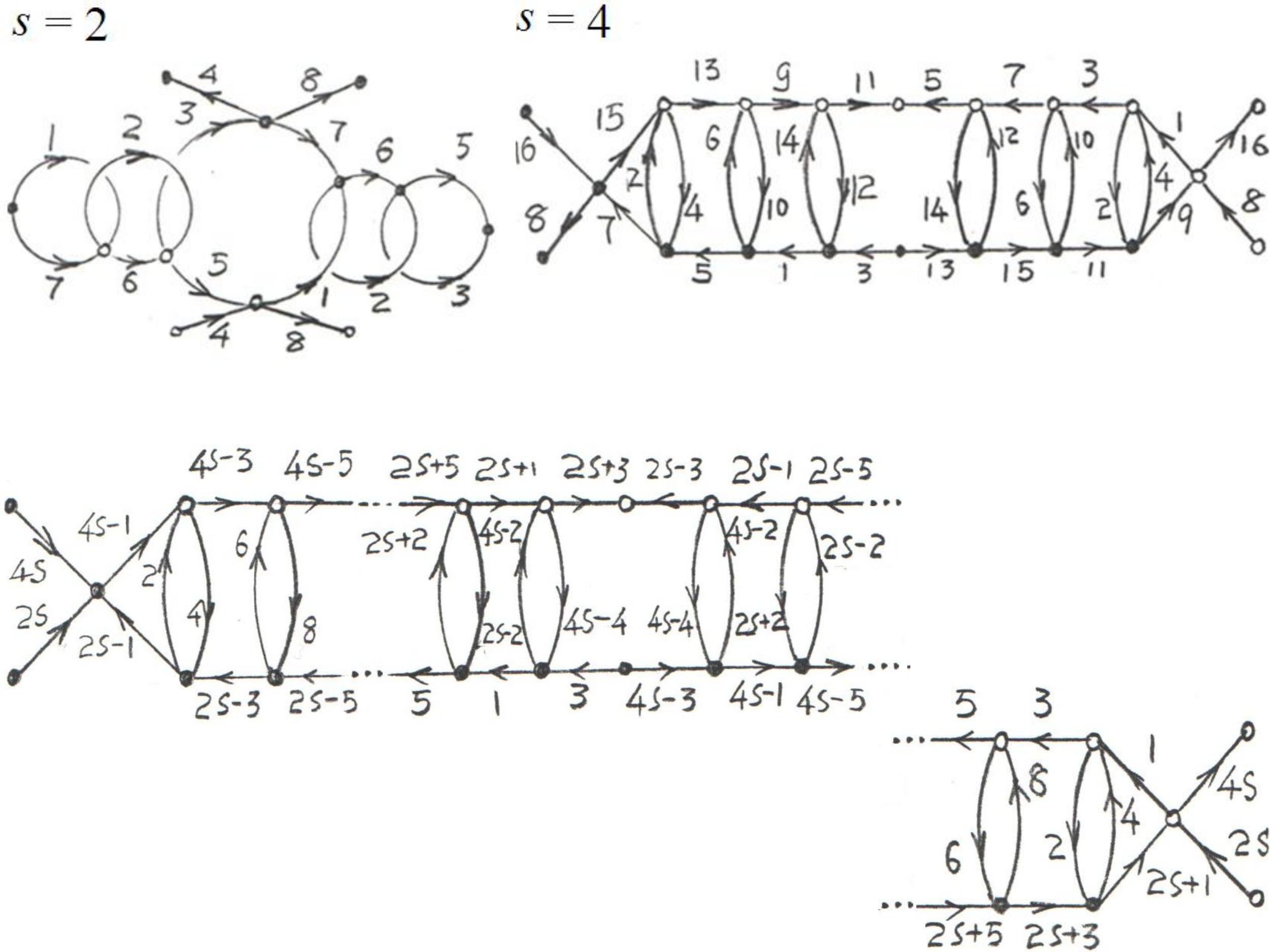}
\end{center}
\begin{center}
    Figure 5.
\end{center}

\bigskip

We now construct current graphs satisfying conditions F1--F7. For $s$  odd, {Figure 4} shows the cases 
$s = 1$,  $ s = 3$,  as well as the case $s$  arbitrary odd. For $s$  even, {Figure 5} shows the cases $s = 2$,   $s = 4$;  it is easy to generalize this current graph for an arbitrary even $s$.

As in the proofs of Lemmas 1 and 2, it can be verified that these index 2 current graphs determine quadrangular embeddings of $K_{8s}$ in $S_{8s^2-5s+1}$.
\end{proof}

\medskip
\noindent
\hbox{\bfseries Lemma 4.}
{\itshape
There exist quadrangular embeddings of $K_{8s+1}-K_4$ in $S_{8s^2-3s-1}$ for $s \geqslant 2$.
}
\begin{proof}
Consider the index 1 current graph for ${\mathbb Z}_{8s-3}$  in Figure 6 with four labeled vertices: $x$, $y$,  $z$,  and $w$.  It is easy to see that for each labeled vertex $v$, the sum of the incident currents is $1$ or $8s-4$, which both are generators of ${\mathbb Z}_{8s-3}$.  Thus, such current graphs determine quadrangular embeddings of $K_{8s+1} - K_{4}$  by Lemmas 1 and 2. 
\end{proof}

\begin{center}
\hspace*{0cm}    \pdfimage width 0.8\textwidth {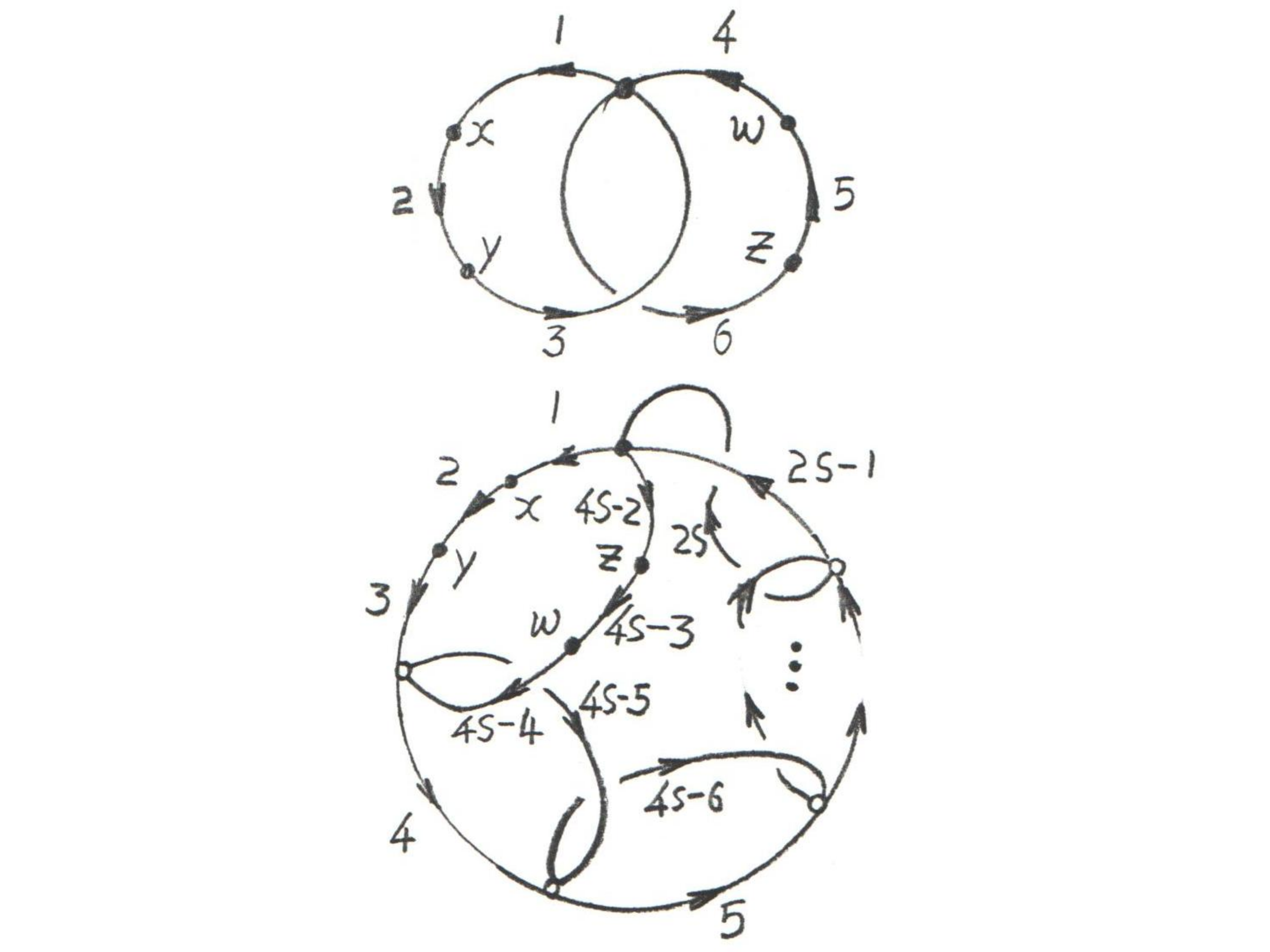}
\end{center}
\begin{center}
    Figure 6.
\end{center}

\begin{center}
\hspace*{0cm}    \pdfimage width 0.71\textwidth {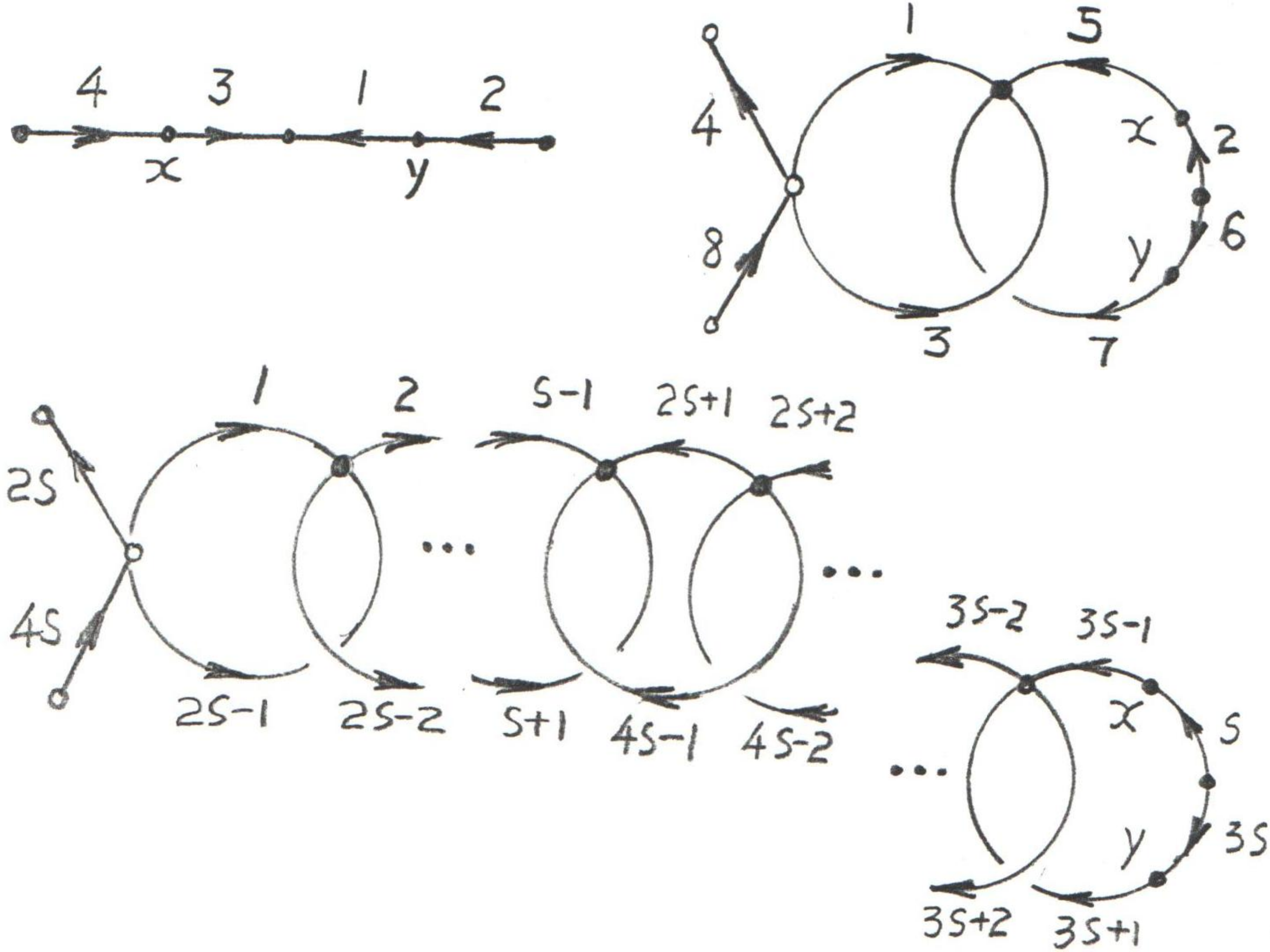}
\begin{center}
\end{center}
    Figure 7.
\end{center}

\noindent
\hbox{\bfseries Lemma 5.}
{\itshape
There exist quadrangular embeddings of  $K_{8s+2}-K_{2}$  in $S_{8s^2-s}$ \, for $s \geqslant 1$.
}
\begin{proof}
For $s = 1$, the index 1 current graph with two labeled vertices $x$  and $y$  in Figure~7 (top) determines a quadrangular embedding of $K_{10}-K_{2}$.
When $s \geqslant 2$, consider the index 1 current graph for the group ${\mathbb Z}_{8s}$ with two labeled vertices $x$  and $y$  in Figure~7 (bottom). For $x$,  the sum of the incident currents is $1$,  which is a generator of ${\mathbb Z}_{8s}$.  For $y$,  the sum of the incident currents is $2s-1$,  a generator of 
${\mathbb Z}_{8s}$, since $8s(s^2-2s+1)+(-4s^2+6s-1)(2s-1) = 1$.  So such current graphs determine quadrangular embeddings of $K_{8s+2}-K_{2}$  by Lemmas 1 and 2.
\end{proof}

\begin{center}
\hspace*{0cm}    \pdfimage width 0.6\textwidth {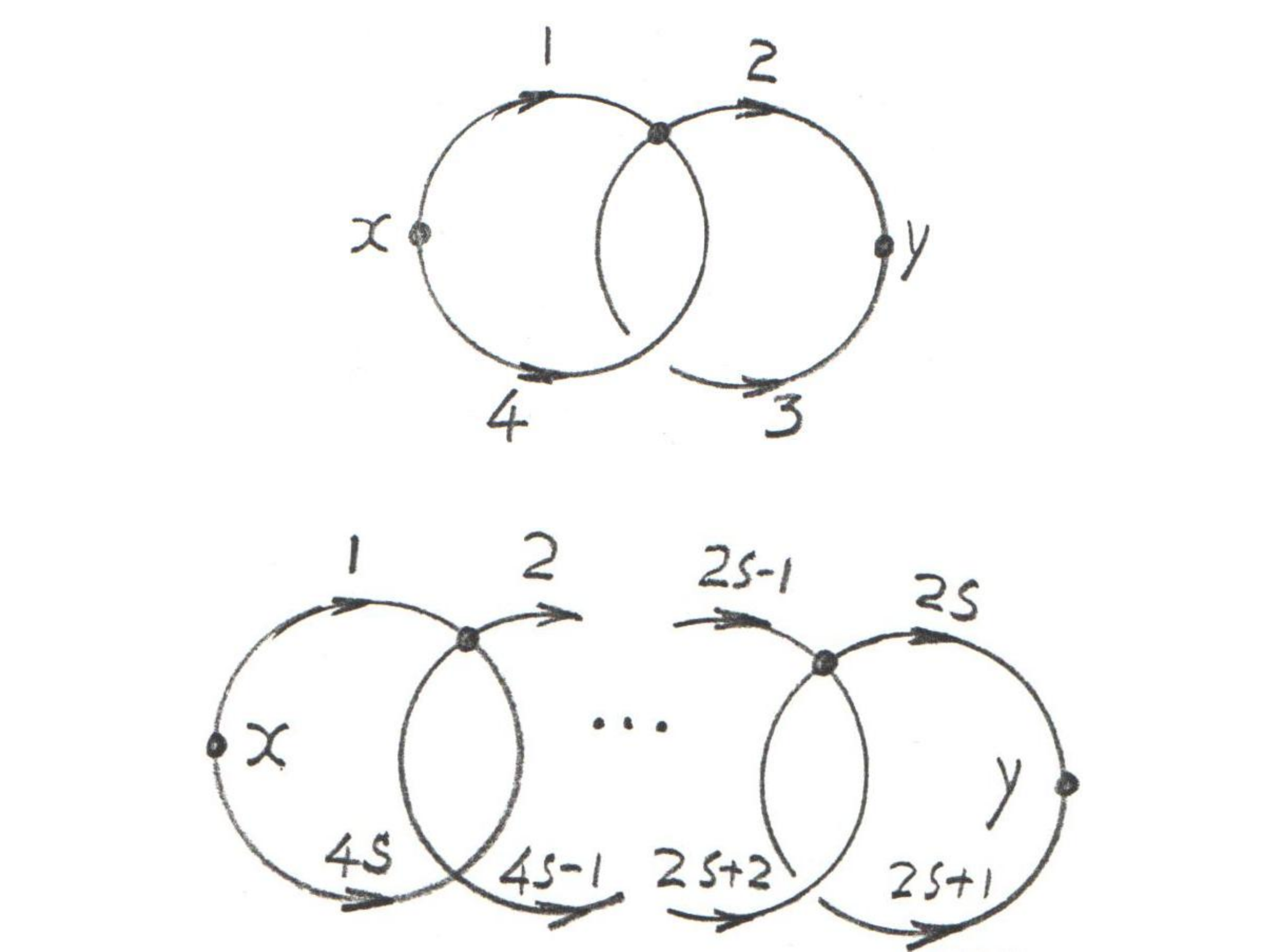}
\end{center}
\begin{center}
    Figure 8.
\end{center}

\medskip

\medskip
\noindent
\hbox{\bfseries Lemma 6}
{\itshape
There exist quadrangular embeddings of $K_{8s+3}-K_{2}$ in $S_{8s^2+s}$ for $s \geqslant 1$.
}
\begin{proof}
Figure 8 shows index 1 current graphs for the group ${\mathbb Z}_{8s+1}$  with two labeled vertices $x$ and $y$. For $x$ and $y$, the sum of the incident currents is $4s + 1$. This is a generator of 
${\mathbb Z}_{8s+1}$ because $-(8s+1)+2(4s+1)=1$. Hence, for $s \geqslant 1$ this graph determines a quadrangular embedding of $K_{8s+3}-K_{2}$ in $S_{8s^2+s}$.
\end{proof}

\medskip
\noindent
\hbox{\bfseries Lemma 7}
{\itshape
There exist quadrangular embeddings of $K_{8s+4}-K_{4}$ in $S_{8s^2+3s-1}$ \, for $s \geqslant 2$.
}
\begin{proof}
For $s \geqslant 2$ even, consider the index 2 current graph for the group ${\mathbb Z}_{8s}$ shown in Figure 9 with four labeled vertices $x$, $y$, $z$, and $w$. For each labeled vertex, the sum of the incident currents is 1 or $8s-1$, which both are generators of ${\mathbb Z}_{8s}$. So this current graph determines a quadrangular embedding of $K_{8s+4}-K_{4}$ in $S_{8s^2+3s-1}$.

For $s$ odd, consider the index 2 current graph for the group ${\mathbb Z}_{8s}$, shown in Figure 10 with four labeled vertices $x$,  $y$,  $z$,  and $w$.  The sum of the incident currents is as follows: $6s +1$  for $x$, $6s-1$ for $y$, $2s-1$ for $z$, and $2s+1$ for $w$, each of which is a generator of ${\mathbb Z}_{8s}$. So this current graph determines a quadrangular embedding of $K_{8s+4}-K_4$ in $S_{8s^2+3s-1}$.
\end{proof}

\begin{center}
\hspace*{0cm}    \pdfimage width 0.88\textwidth {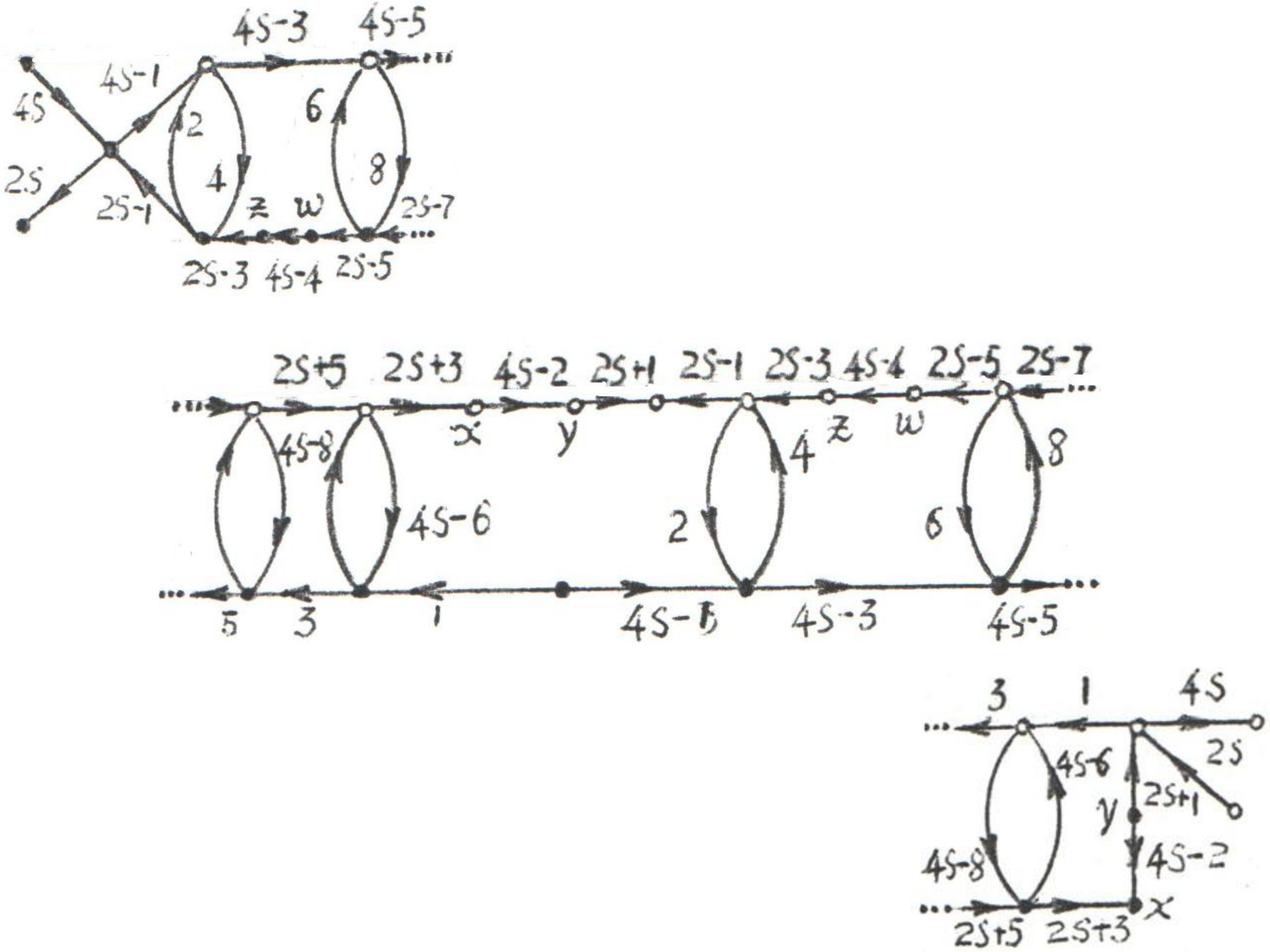}
\end{center}
\begin{center}
    Figure 9.
\end{center}

\smallskip

\begin{center}
\hspace*{0cm}    \pdfimage width 0.88\textwidth {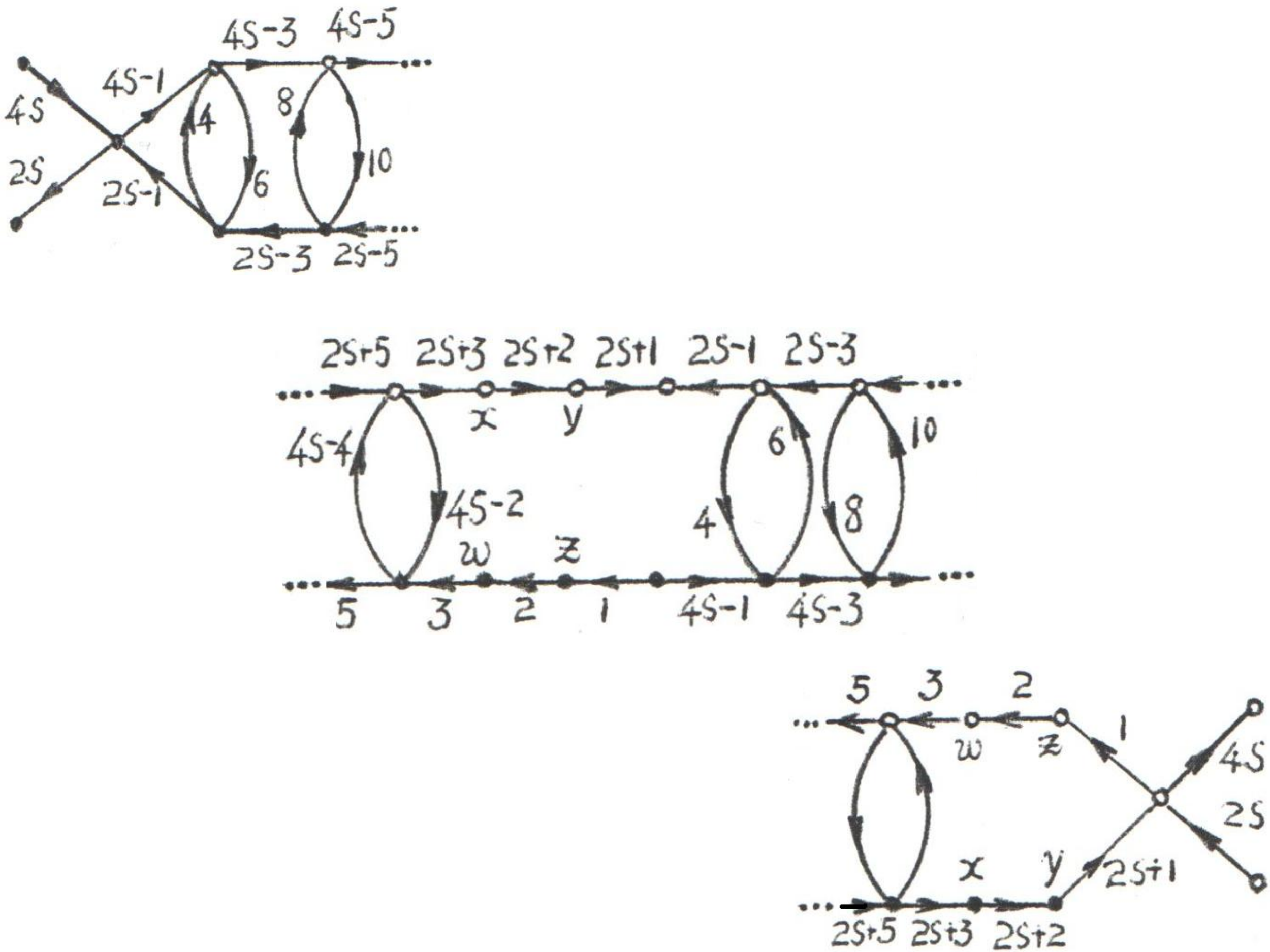}
\end{center}
\begin{center}
    Figure 10.
\end{center}

\bigskip

The following result is due to Archdeacon [1].

\medskip
\noindent
\hbox{\bfseries Lemma 8}
{\itshape
There exist quadrangular embeddings of $K_{8s+5}$ in $S_{8s^2+5s+1}$ \, for $s \geqslant 0 $.
}

\medskip

\begin{center}
\hspace*{0cm}    \pdfimage width 0.73\textwidth {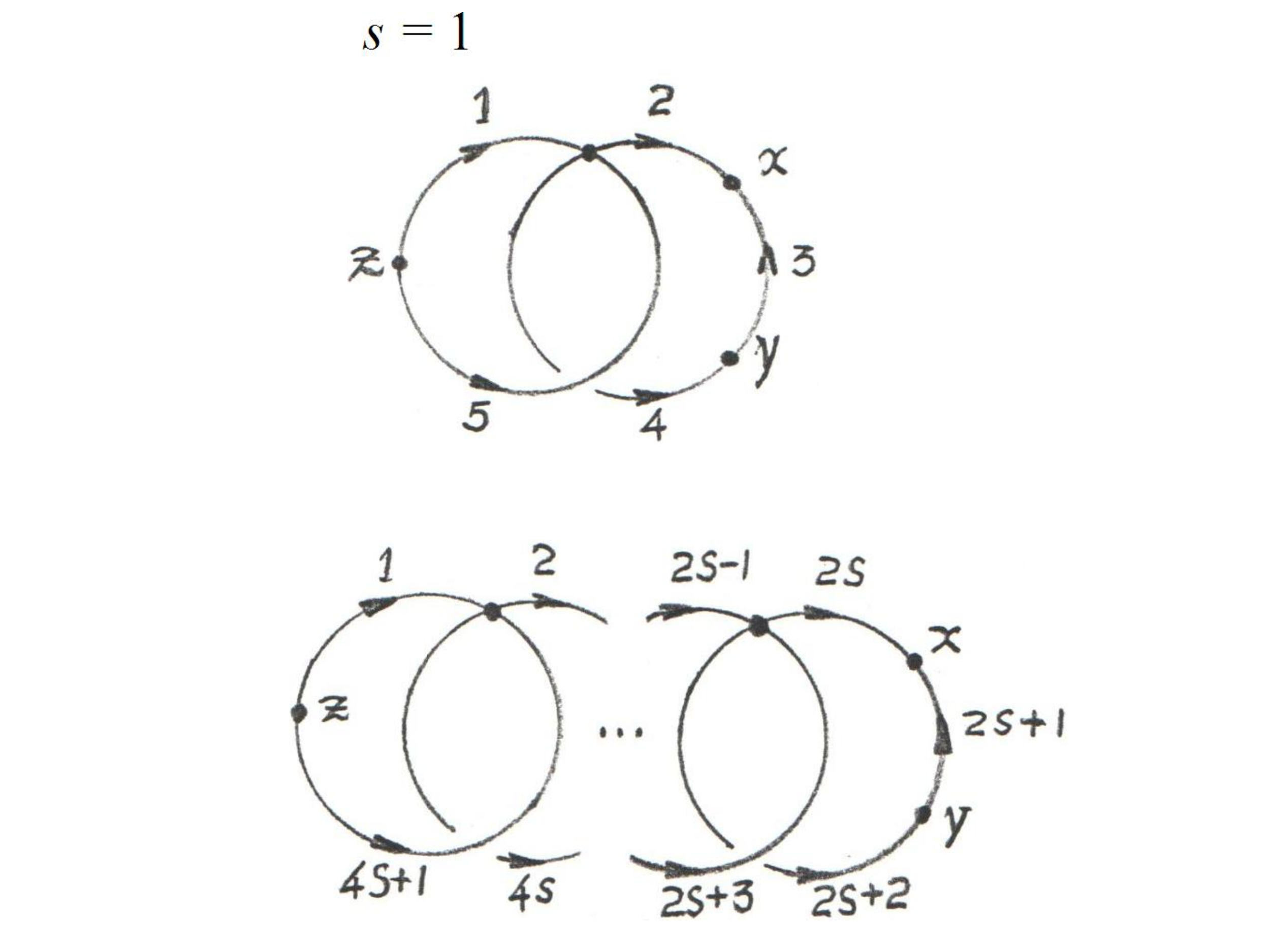}
\end{center}
\begin{center}
    Figure 11.
\end{center}

\bigskip

\begin{center}
\hspace*{0cm}    \pdfimage width 0.73\textwidth {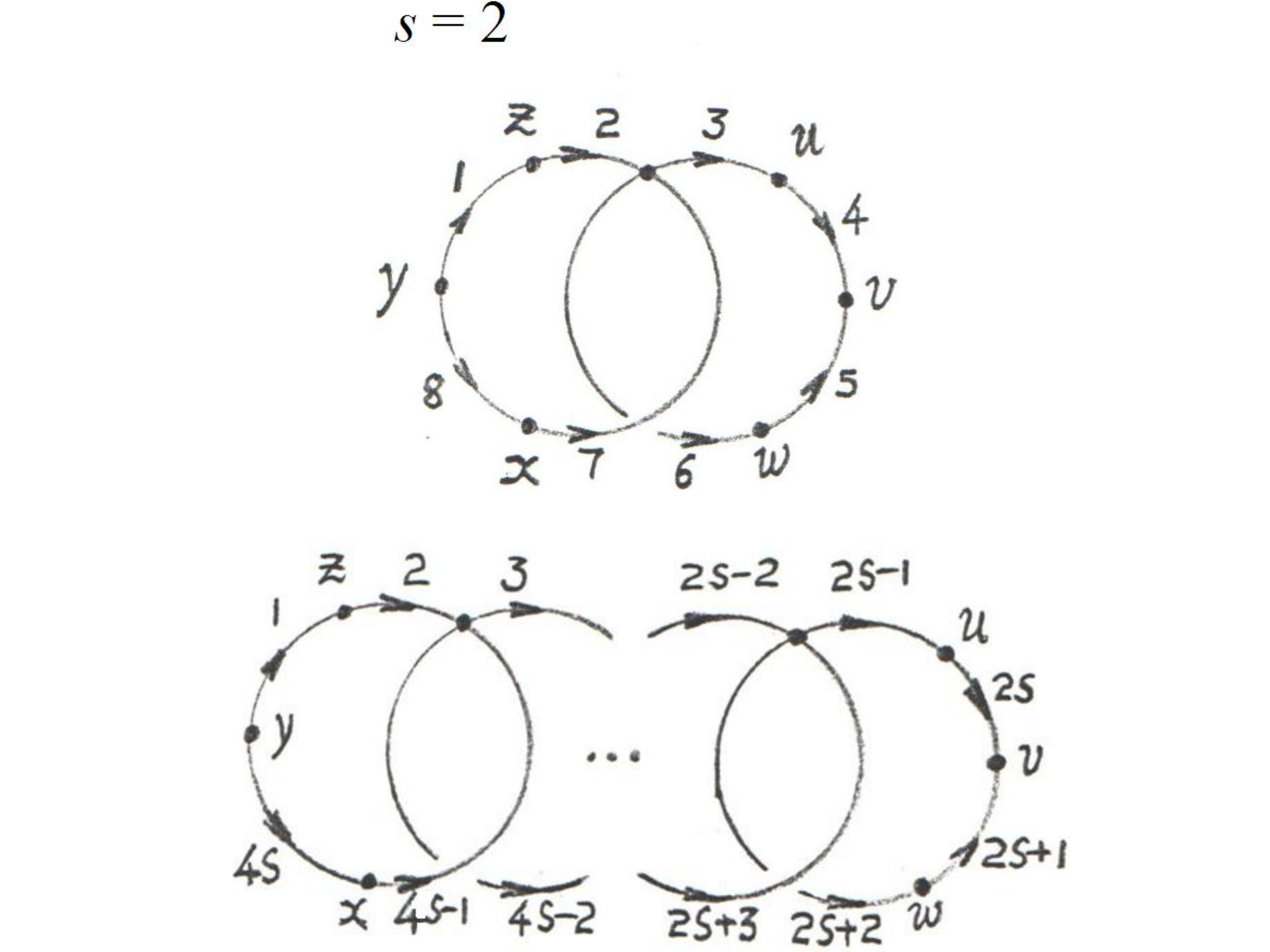}
\end{center}
\begin{center}
    Figure 12.
\end{center}

\bigskip

\noindent
\hbox{\bfseries Lemma 9}
{\itshape
There exist quadrangular embeddings of $K_{8s+6}-K_{3}$ in $S_{8s^2+7s+1}$ \, for $s \geqslant 1$.
}
\begin{proof}
Consider the index 1 current graph for the group ${\mathbb Z}_{8s+3}$ with three labeled vertices $x$, $y$, and $z$ depicted in Figure 11. For both $x$ and $y$, the sum of the incident currents is $4s+1$, which is a generator of ${\mathbb Z}_{8s+3}$. For $z$, the sum of the incident currents is 1, which also is a generator. So this current graph determines a quadrangular embedding of $K_{8s+6}-K_{3}$ in $S_{8s^2+7s+1}$.
\end{proof}

\medskip
\noindent
\hbox{\bfseries Lemma 10}
{\itshape
There exist quadrangular embeddings of $K_{8s+7}-K_{6}$ in $S_{8s^2+9s-1}$ \, for $s \geqslant 2$.
}
\begin{proof}
Consider the index 1 current graph for the group ${\mathbb Z}_{8s+1}$ with six labeled vertices $x$, $y$, $z$, $u$, $v$, and $w$ in Figure 12. For each of $x$, $z$, $u$, and $w$, the sum of the incident currents is 1 or $8s$ which both are generators of ${\mathbb Z}_{8s+1}$.  For $y$  and $v$,  the sum of the incident currents is $4s$  or $4s+1$.  Thus, these current graphs determine quadrangular embeddings of $K_{8s+7}-K_{6}$  in $S_{8s^2+9s-1}$. 
\end{proof}

\medskip
\section{Proof of Theorem 1}

Given an embedding of a graph in a surface, denote by $r_i$  the number of faces of length $i$,  where $i = 3, 4, 5, \ldots$. 
An $(r_3, r_4, r_5, \ldots)${\it -embedding} is a 2-cell embedding having exactly $r_i$  faces of length 
$i$.

First, we show that $\gamma_4 (K_p) \geqslant \lceil p(p - 5)/8 \rceil + 1$.  By assumption, we have an 
$(r_4, r_5, r_6, \ldots)$-embedding of $K_p$  in $S_{\gamma_4}(K_p)$.  Denote by $p$, $q$, and $r$ the numbers of vertices, edges, and faces in that embedding. By Euler's formula, we have:
\begin{equation}
p - q + r = 2 - 2 \gamma_4 (K_p).
\end{equation}
\noindent Furthermore, we have:
\begin{equation}
r_4 + r_5 + r_6 + \cdots = r
\end{equation}
and 
\begin{equation}
4 r_4 + 5 r_5 + 6 r_6 = 2q.
\end{equation}
Hence, 
\begin{equation}
\sum_{i \geqslant 5} (i - 4) r_i = 2q - 4r.
\end{equation}
\noindent From Eqs. (1) and (4) we deduce
$$
\gamma_4 (K_p) = \frac {p(p - 5)} 8 + 1 + \frac 1 8 \sum _{i \geqslant 5} (i  - 4 )r_i \geqslant  \Big \lceil  \frac {p(p - 5)} 8 \Big \rceil + 1.
$$

To complete the proof, we need to show that $\gamma_4 (K_p) \leqslant \lceil p(p - 5)/8 \rceil + 1$.  This task splits into a number of cases.

\begin{center}
\hspace*{0cm}    \pdfimage width 1.0\textwidth {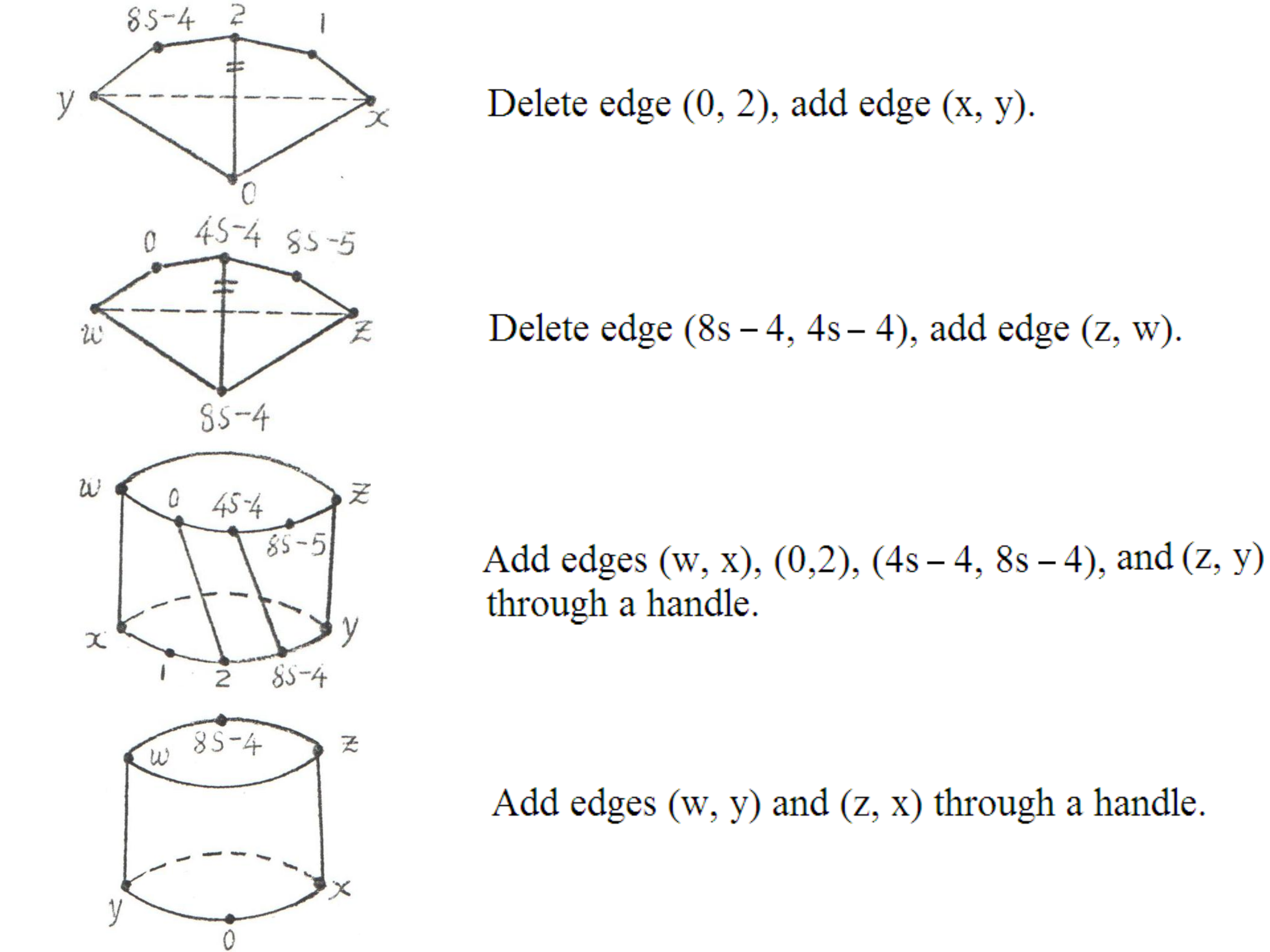}
\end{center}
\begin{center}
    Figure 13.
\end{center}

\bigskip

\begin{center}
\hspace*{0cm}    \pdfimage width 0.7\textwidth {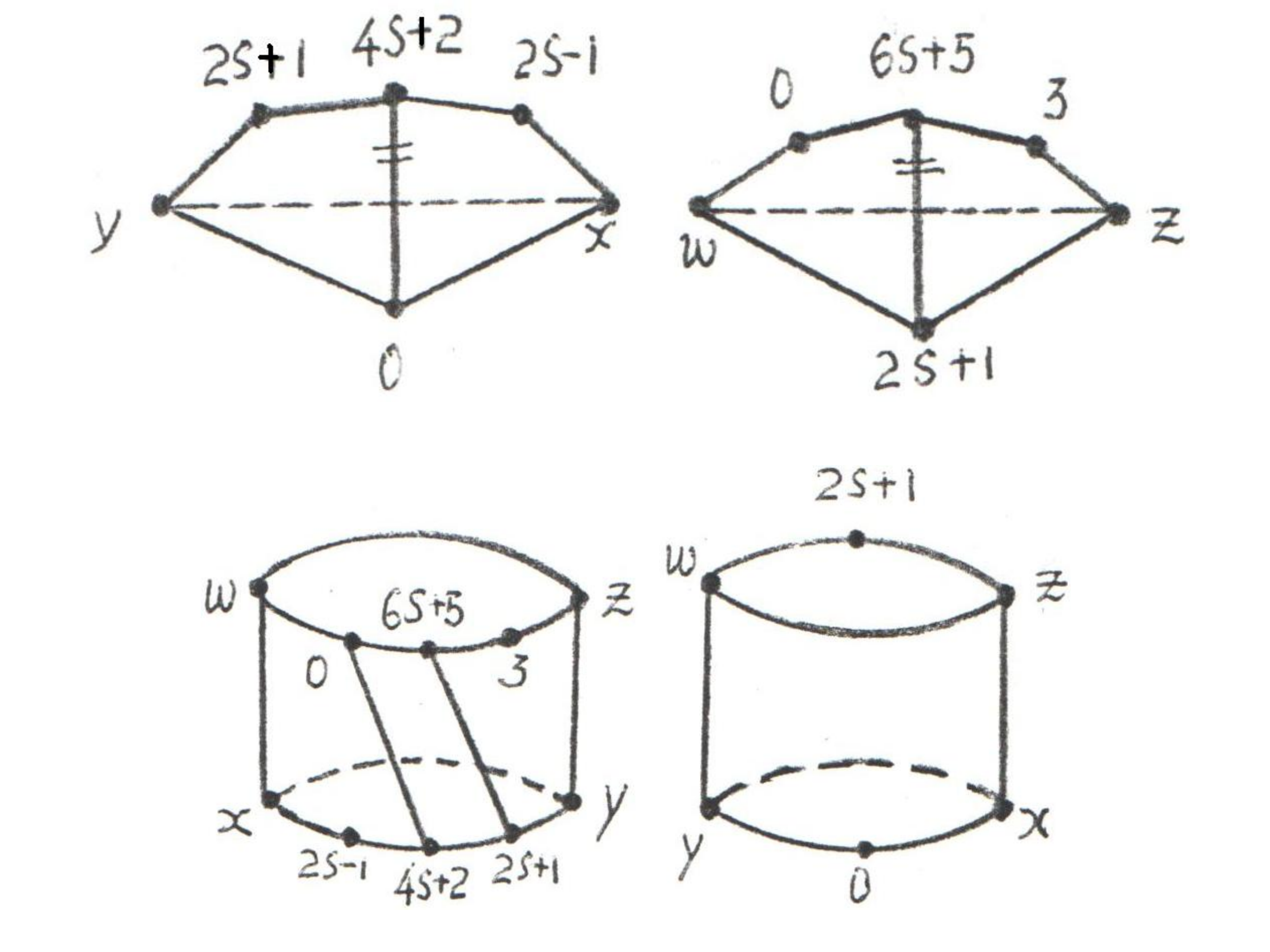}
\end{center}
\begin{center}
    Figure 14.
\end{center}

\begin{center}
\hspace*{0cm}    \pdfimage width 0.65\textwidth {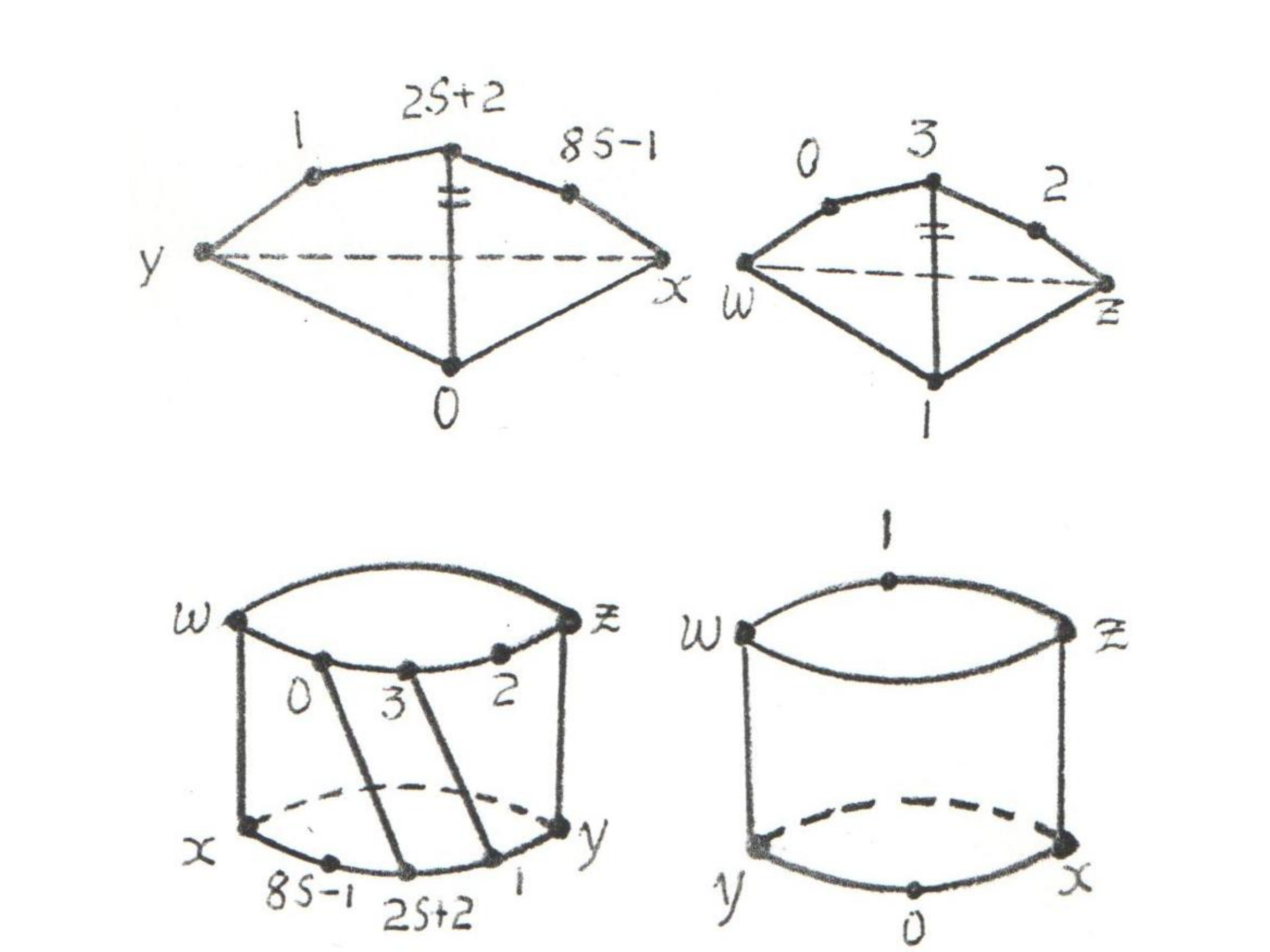}
\end{center}
\begin{center}
    Figure 15.
\end{center}

\bigskip

For $p \equiv 0 \, \, \, \, {\rm {or}} \, \, \, \,  5 \mod 8$, by Lemmas 3 and 8, there exist quadrangular embeddings of $K_p$ in $S_{p(p-5)/8 + 1}$.  Hence, $\gamma_4 (K_p) \leqslant p(p - 5)/8 + 1$  for 
$p \equiv 0 \, \, \, \, {\rm {or}} \, \, \, \,  5 \mod 8$.

For $p \equiv 2 \, \, \, \, {\rm {or}} \, \, \, \,  3 \mod 8$ ($p \geqslant 10$), by Lemmas 5 and 6, there exist quadrangular embeddings of $K_p - K_2$ in $S_{\lceil p(p-5)/8 \rceil}$,  and we can add one missing edge with addition of one extra handle to obtain a 2-cell embedding of $K_p$  in $S_{\lceil p(p-5)/8 \rceil + 1}$  with $r_3 = 0$.  Hence, $\gamma_4(K_p) \leqslant \lceil p(p-5)/8 \rceil +1$.

For $p \equiv 1 \, \, \, \, {\rm {or}} \, \, \, \,  4 \mod 8$ \, ($p \ne 9 \, \, \, \, {\rm {or}} \, \, \, \,  12$),  by Lemmas 4 and 7, there exist quadrangular embeddings of $K_p - K_4$  in 
$S_ {\lceil p(p-5)/8 \rceil - 1}$ .  All we need is to add six missing edges with addition of two extra handles to have $K_p$  embedded in $S_{\lceil p(p-5)/8 \rceil +1}$  with $r_3 = 0$.  Here we have to consider three cases as well as the exceptional graphs  $K_9$  and  $K_{12}$.

\medskip
\noindent {\it Case 1. $p = 8s + 1$.} From the construction of quadrangular embeddings in the proof of Lemma 4, we have a face distribution as in Figure 13. Then we add missing edges with addition of two extra handles as in Figure 13.

\medskip
\noindent {\it Case 2. $p = 8s + 4$, where $s$ is odd not less than $3$.} From the construction of quadrangular embeddings in Lemma 5, we have a face distribution as in Figure 14 which indicates the procedure of addition of missing edges.

\medskip
\noindent {\it Case 3. $p = 8s + 4$, where $s$ is even.} We have a face distribution as in Figure 15 which indicates the procedure of addition of missing edges in this case.
\medskip

For $K_9$,  we first construct an embedding of $K_9 - K_3$ in $S_4$ by using an index 1 current graph for the group ${\mathbb Z}_6$, with three labeled vertices $x$, $y$, and $z$. Row $0$ is as follows:

$$
0 : x, 1, y, 2, z, 3, 4, 5.
$$

\noindent This embedding has six faces of length 3. We add three missing edges with addition of two extra handles as in Figure 16.
\bigskip

\begin{center}
\hspace*{0cm}    \pdfimage width 0.7\textwidth {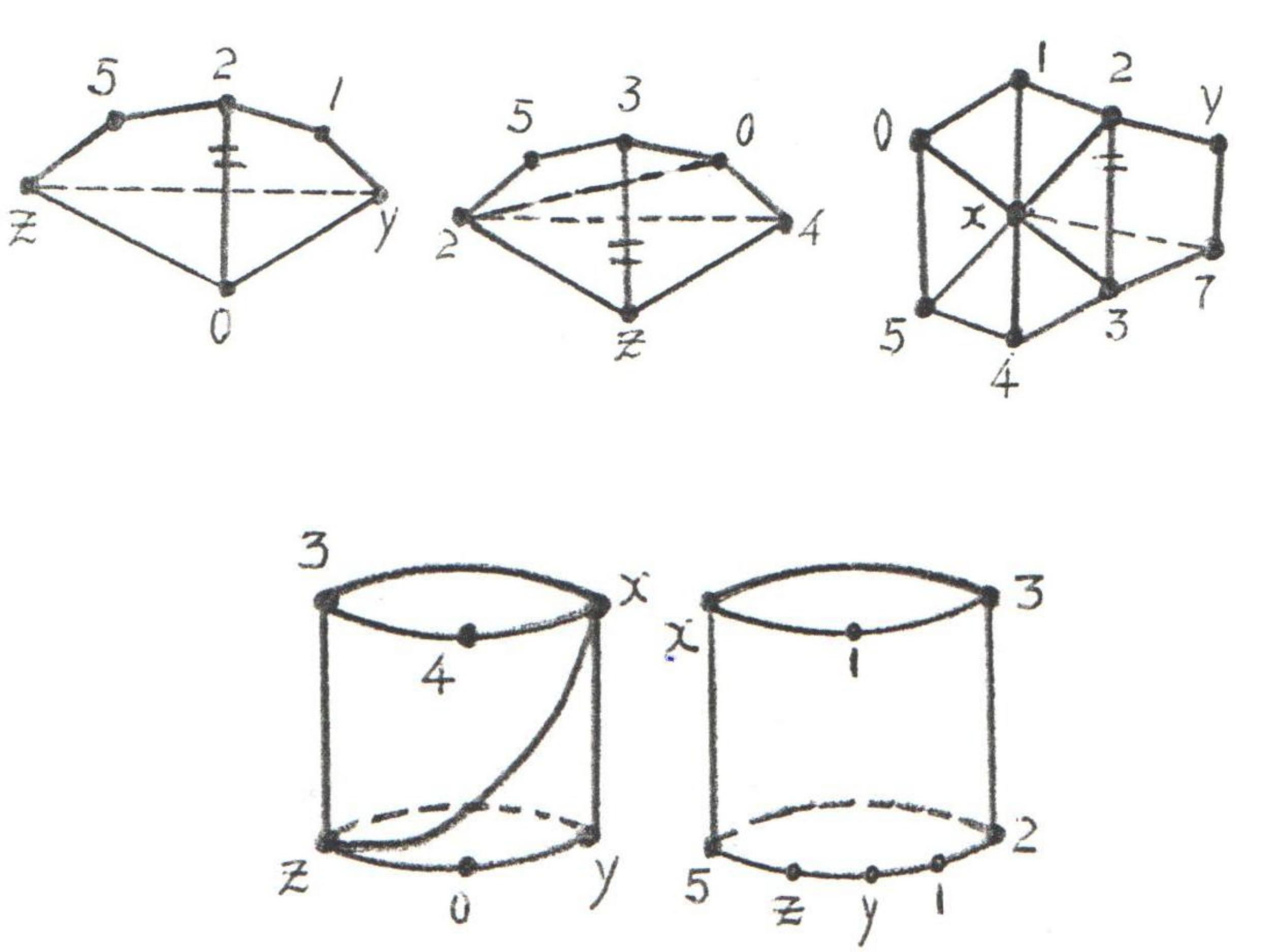}
\begin{center}
\end{center}
    Figure 16.
\end{center}

\begin{center}
\hspace*{0cm}    \pdfimage width 0.8\textwidth {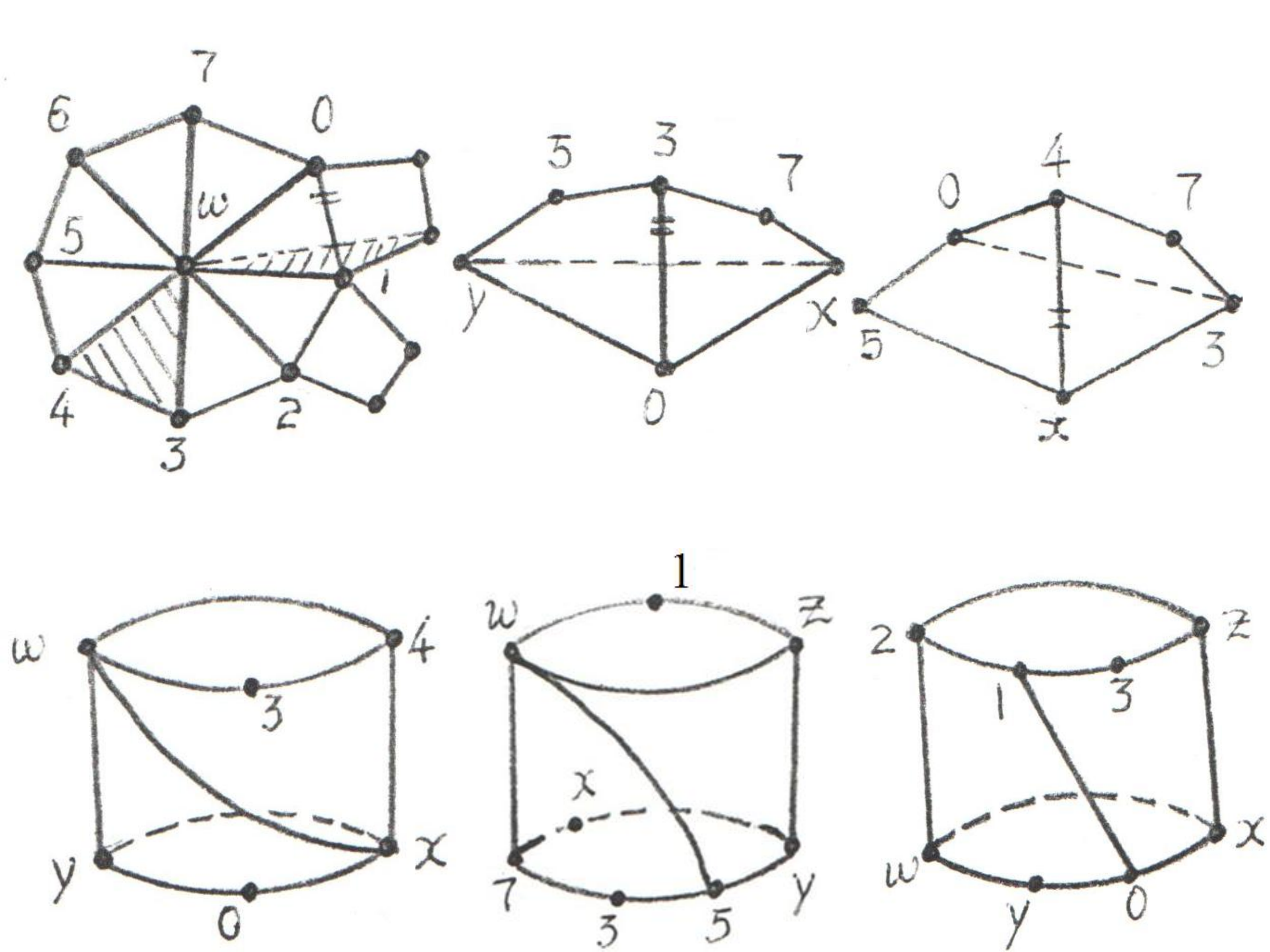}
\begin{center}
\end{center}
    Figure 17.
\end{center}

\begin{center}
\hspace*{0cm}    \pdfimage width 0.7\textwidth {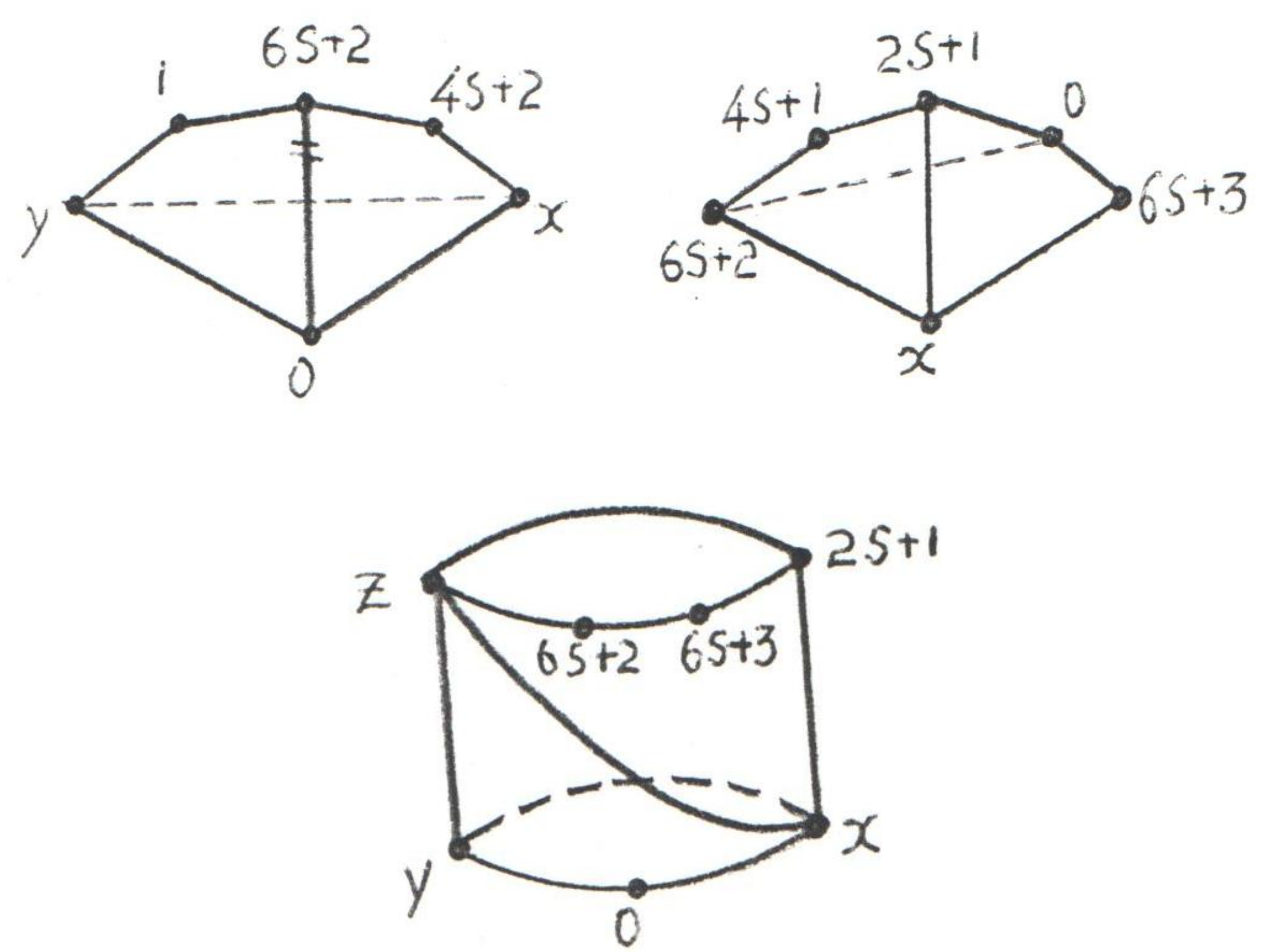}
\end{center}
\begin{center}
    Figure 18.
\end{center}

\bigskip

\begin{center}
\hspace*{0cm}    
\pdfimage width 0.7\textwidth {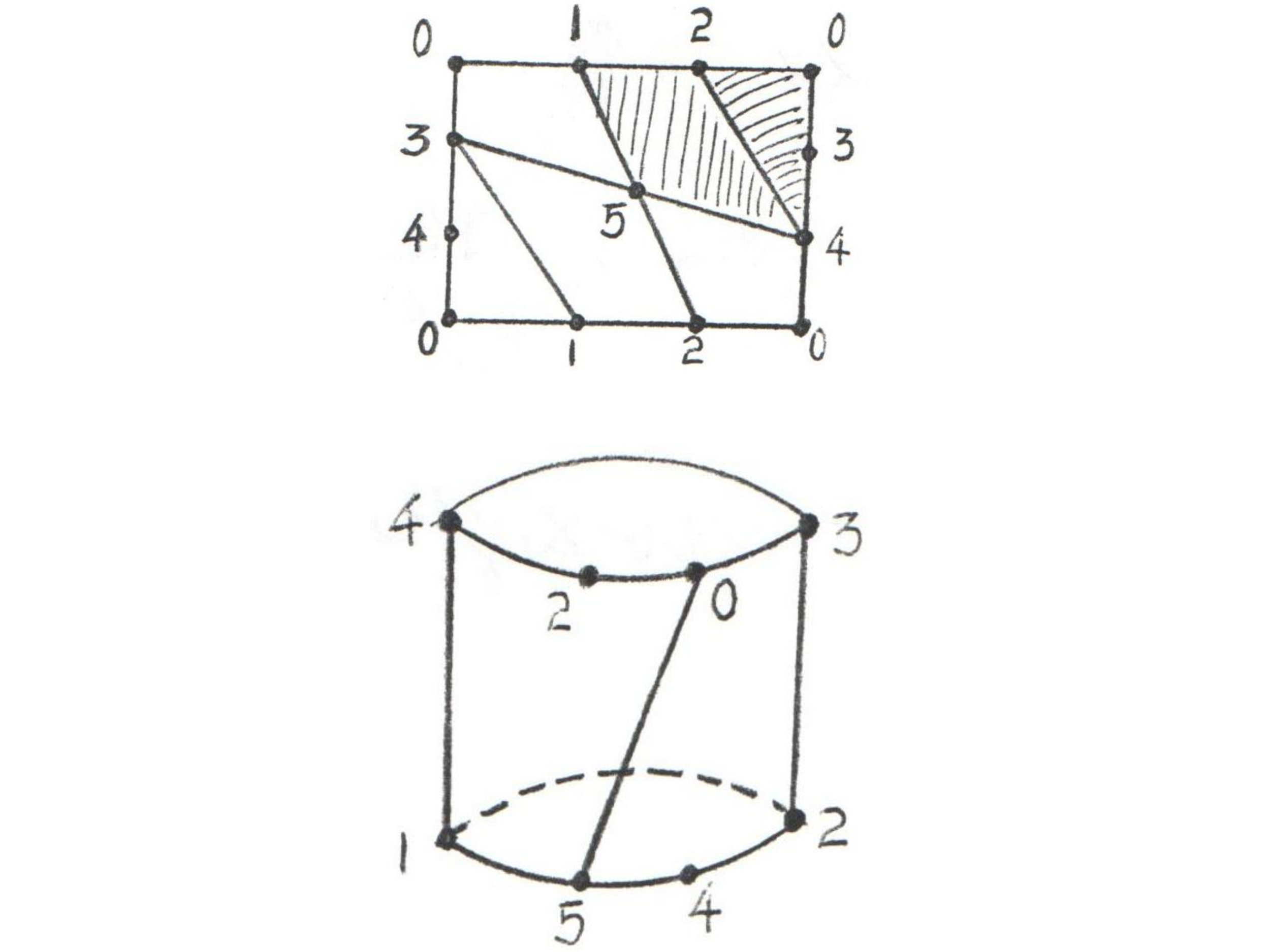}
\end{center}
\begin{center}
    Figure 19.
\end{center}

For $K_{12}$,  we first construct an embedding of $K_{12} - K_4$  in $S_9$  with $r_3 = 8$,  $r_4 = 24$  by using an index 1 current graph for $\mathbb{Z}_8$  with four labeled vertices. Row $0$  is as follows:

$$
0 : 4, x, 3, y, 2, z, 1, w, 7, 6, 5.
$$

\noindent Then we add six missing edges with addition of three extra handles as in Figure 17.

\begin{center}
\hspace*{0cm}    
\pdfimage width 0.7\textwidth {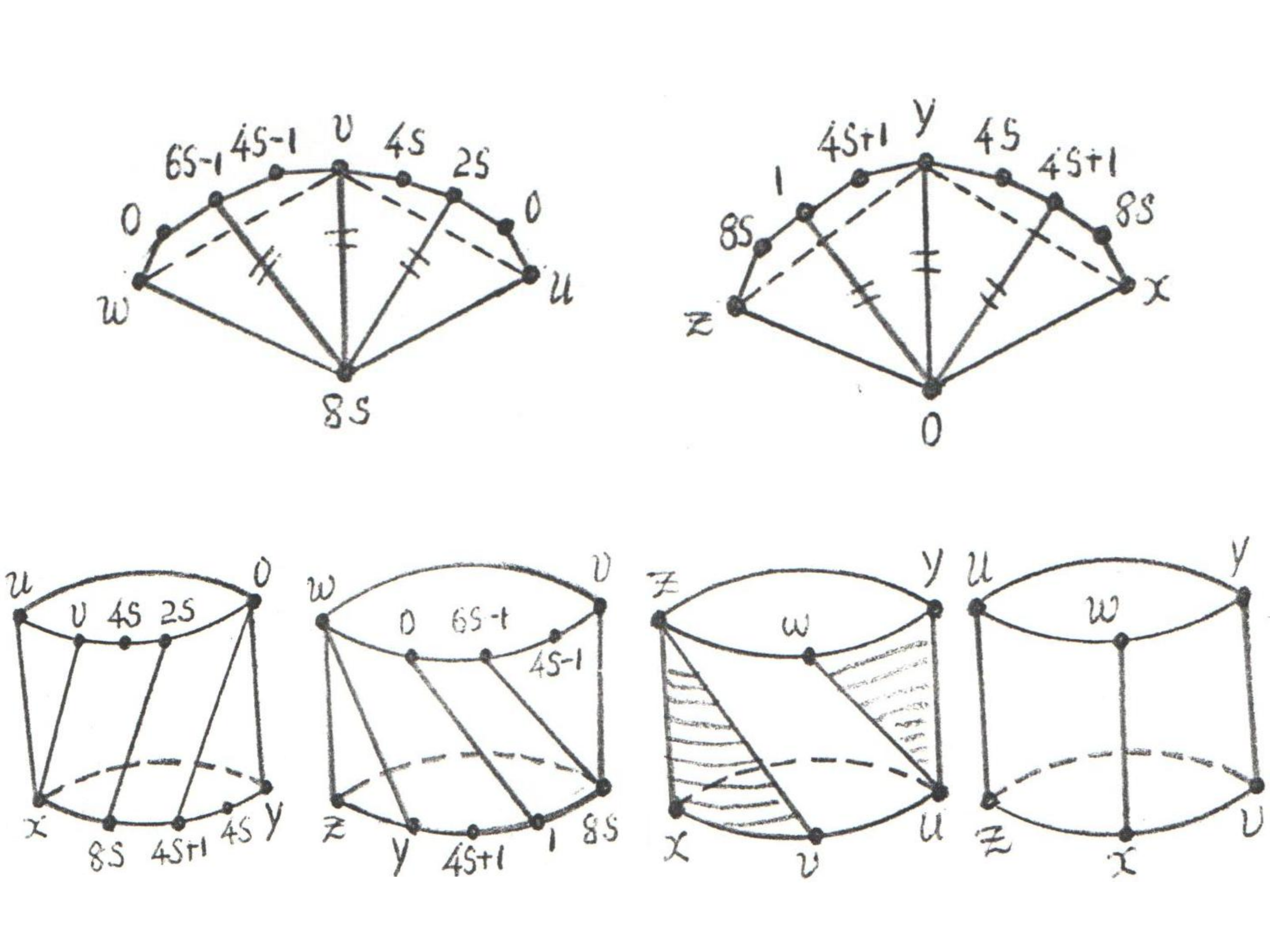}
\end{center}
\begin{center}
    Figure 20.
\end{center}

\bigskip

\begin{center}
\hspace*{0cm}    
\pdfimage width 0.7\textwidth {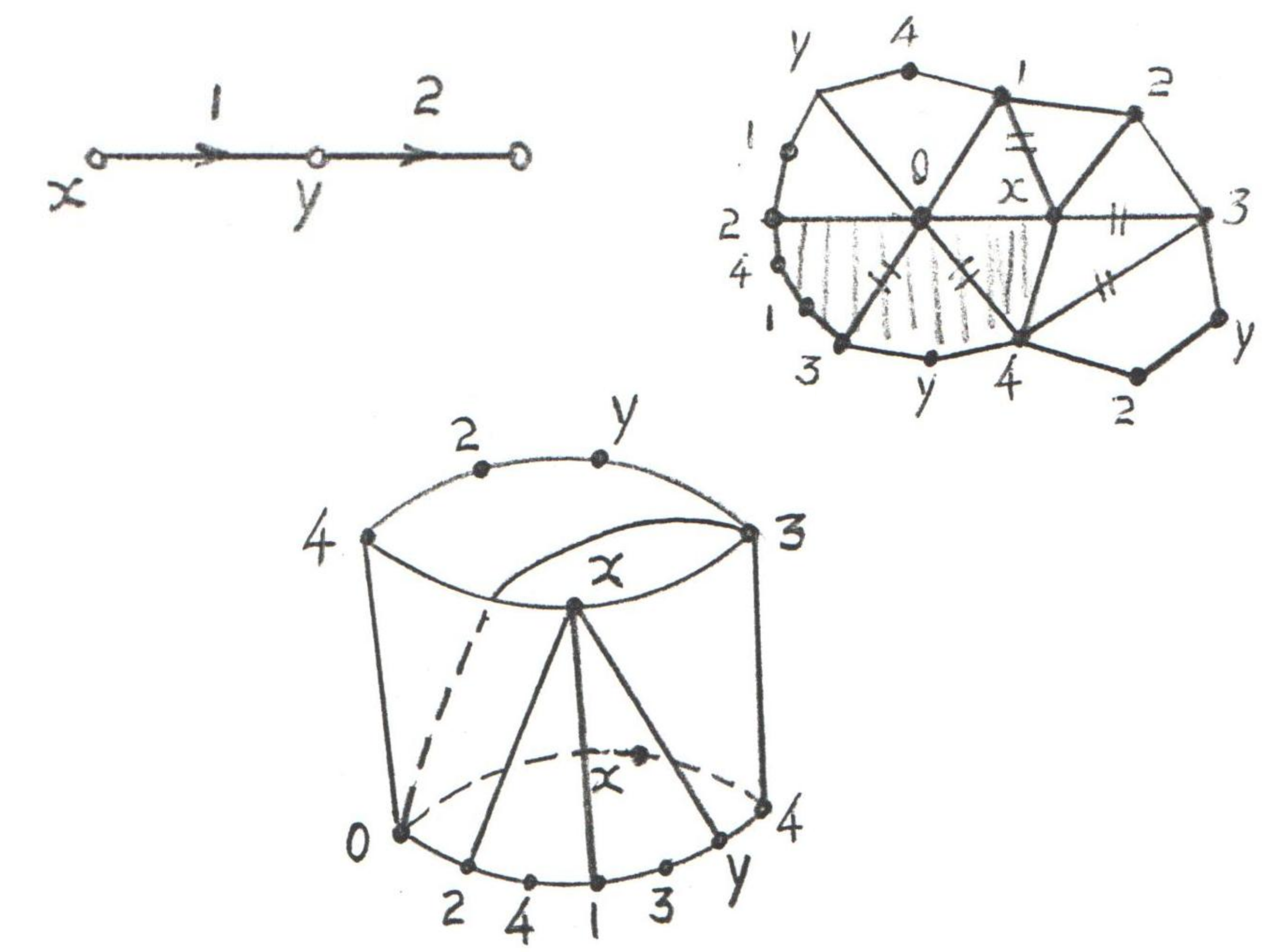}
\end{center}
\begin{center}
    Figure 21.
\end{center}

\medskip

For $p \equiv 6 \mod 8$ ($p \ne 6$), by Lemma 9, there exists a quadrangular embedding of $K_p - K_3$  in 
$S_{\lceil p(p-5)/8 \rceil}$.  This quadrangulartion has a face distrbution as in Figure 18, which shows how to add three missing edges with addition of one extra handle. For $K_6$,  we first embed 
$K_6 - K_3$  in $S_1$  as in Figure 19 (top). Then we add three missing edges with addition of one extra handle as in Figure 19 (bottom).

For $p \equiv 7 \mod 8$ \, ($p \ne 7$),  by Lemma 10, we have quadrangular embeddings of $K_p - K_6$ in
$S_{\lceil p(p-5)/8 \rceil - 3}$. We add fifteen missing edges with additin of four extra handles as in Figure 20. For $K_7$,  we embed $K_7 - K_2$ in $S_2$  by using the index 1 current graph shown in Figure 21. In this embedding, $r_3 = 5$, $r_4 = 4$, and $r_5 = 1$. The face distribution is as in Figure 21 which also shows how we add one missing edge with addition of one handle. 

The proof of Theorem 1 is complete.

\section*{References}

\noindent 1. D. Archdeacon, The medial graph and voltage-current duality, Discrete Math. 104 (1992) 111--141.

\noindent 2. N. Hartsfield and G. Ringel, Minimal quadrangulation of orientable surface, J. Combin. Theory Ser. B 46 (1989) 84--93.

\noindent 3. G. Ringel and J. W. T. Youngs, Solution of the Heawood map-coloring problem, Proc. Nat. Acad. Sci., USA 60 (1968) 438--445.

\noindent 4. A. T. White, Topological graph theory, in: ``Selected Topics in Graph Theory,'' Eds. L. W. Beineke and R. J. Wilson, Academic Press, London, 1978, pp. 15--50.

\noindent 5. A. T. White, The Proof of Heawood conjecture, in: ``Selected Topics in Graph Theory,'' Eds. L. W. Beineke and R. J. Wilson, Academic Press, London, 1978, pp. 51--82.

\medskip
\section{Appendix}

The genus of the complete graph was established by Ringel and Youngs \cite{R} and was mainly concerned with triangulations of surfaces. Nonetheless, since then a great deal of interest has also been generated in quadrangulations of surfaces. In particular, Hartsfield and Ringel \cite{HR, HR2}, the first author (Lawrencenko) \cite{L2, L3} and Lawrencenko et al. in the original 1998 version of the current paper have considered minimal quadrangulations of surfaces. The purpose of this appendix is to identify and clarify copyright issues around the quadrangular genus of complete graphs and related topics. 

Similarly to the orientable case, we define the nonorientable quadrangular genus of $G$, denoted $\gamma^{*}_4 (G)$, to be the minimum value of $k$ for which $G$ has a 2-cell embedding in $N_k$ (the sphere with $k$ crosscaps) such that the smallest face in the embedding is a quadrangle. 
In Sections 1--4, we established the quadrangular genus for the complete graph for the orientable case
(Theorem 1) while Nora Hartsfield did the nonorientable case:

\medskip
\noindent
\hbox{\bfseries Theorem 2 (Hartsfield).}
{\itshape
$$
\gamma^{*}_4 (K_p) = \bigg{\lceil} \frac {p(p-5)}{4} \bigg{\rceil} +2.
$$
}
\medskip

In fact, the original 1998 version of the current paper was submitted to Discrete Mathematics in June 1998 but shortly after that the first author (Lawrencenko) withdrew it from the journal. The reason for the withdrawal was that the first author learned from Nora Hartsfield that she had obtained Theorem 2. 
Then, Nora Hartsfield and the first author agreed to merge their papers and produce one joint paper with both Theorems 1 and 2 instead of two complementary papers drafted separately. Unfortunately, the first author was too slow in fulfilling the project, and as a result no journal paper has been published as yet. Sadly, Nora passed away in 2011. The first author learned about her death only in winter 2016.

It must be emphasized that the main result of paper \cite{LEYZ} authored by Liu~W., Ellingham~M.\,N., Ye~D., and Zha~X. (the first author's former co-workers) immediately follows from Theorems 1, 2 of the current paper.
In addition, the original 1998 version of the current paper has been circulating over the internet around the mathematics community; for instance, it was used by Vladimir Korzhik, Yusuke Suzuki, etc.; some of the citations of it are present in the reference list of Suzuki's important paper \cite{S}.
Moreover, Theorem 1 has been used by Yusuke Suzuki \cite{S}, with proper referencing to the original 1998 version of the current paper, for Suzuki's proofs of the following two of his results. The real significance of these results is that they perfectly show the relationship between triangulations and quadrangulations of closed surfaces. 

\medskip
\noindent
\hbox{\bfseries Theorem 3 (Suzuki \cite{S}).}
{\itshape There exists an integer $h_0$ such that for any two surfaces $S_{h_1}$ and $S_{h_2}$ satisfying the inequalities $h_1 \geqslant h_0$ and 
\medskip

$$
2h_2 - 3h_1 -\bigg{\lceil} \frac {-331 + 19 \sqrt{1+48h_1}}{12} \bigg{\rceil} \geqslant -1,
$$
\noindent there exists a triangulation of $S_{h_1}$ whose graph quadrangulates $S_{h_2}$}.
\medskip

\noindent
\hbox{\bfseries Theorem 4 (Suzuki \cite{S}).}
{\itshape There exists an integer $k_0$ such that for any two surfaces $N_{k_1}$ and $S_{h_2}$ satisfying the inequalities $k_1 \geqslant k_0$ and 
\medskip

$$
2h_2 - k_1 -\bigg{\lceil} \frac {-331 + 19 \sqrt{1+24 k_1}}{12} \bigg{\rceil} \geqslant -1,
$$
\noindent there exists a triangulation of $N_{k_1}$ whose graph quadrangulates $S_{h_2}$}.

\medskip
However, Suzuki \cite{S} mistakenly attributes Theorem 2 to Hartsfield and Ringel \cite{HR}: Paper \cite{HR} covers only the case $n \equiv 1 \mod 4$ and the proof of the following theorem should be built upon Theorem 2.

\medskip
\noindent
\hbox{\bfseries Theorem 5 (Suzuki \cite{S}).}
{\itshape There exists an integer $h_0$ such that for any two surfaces $S_{h_1}$ and $N_{k_2}$ satisfying the inequalities $h_1 \geqslant h_0$ and 
\medskip

$$
k_2 - 3h_1 - \bigg{\lceil} \frac {-91 + 11 \sqrt{1+48 h_1}}{12} \bigg{\rceil} \geqslant -1,
$$
\noindent there exists a triangulation of $S_{h_1}$ whose graph quadrangulates $N_{k_2}$}.
\medskip

The question arises about the geometric realizability of quadrangulations of closed surfaces in Euclidean space $\mathbb{R}^d$. Notice that the boundary 2-complex of the 3-cube on 8 vertices can be realized as a regular hexahedron. Furthermore, the regular triangulation of the torus on 8 vertices (which is a subcomplex of the 2-skeleton of a 16-cell) is geometrically realized in $\mathbb{R}^3$ \cite{L}; and, moreover, by removing a suitable set of eight edges from the realization, one can obtain a ``quasi-geometric'' realization in the form of a toroidal polyhedron with all faces quadrilaterals possibly bent about their diagonal axes. The least likely candidates for a positive answer to the quasi-geometric realizability question are the quadrangulations by complete graphs, even with allowing for not only diagonal bends but also curved quadrilateral faces (but with all bounding edges straight). More specifically, {\it is any quadrangulation of $S_4$ with the complete graph on 8 vertices quasi-geometrically realizable in $\mathbb{R}^3$  (or  $\mathbb{R}^d$)?}
\medskip

\noindent {\bf Acknowledgments.} The first author (Lawrencenko) is grateful to Yusuke Suzuki and Vladimir Korzhik for their interest in this work.

\par\medskip 
\par\medskip 
\par\medskip 

\noindent \textsc{S. Lawrencenko} \\
\medskip
\noindent E-mail: \url{lawrencenko@hotmail.com} \\
\noindent Current affiliation: Russian State University of Tourism and Service,\\
\noindent Institute for Tourism and Hospitality,\\
\noindent Bldg. 32A, Kronstadt Boulevard, Moscow, 125438, Russia
\smallskip

\noindent (Former Affiliation: Department of Mathematics, \\
\noindent Hong Kong University of Science and Technology, \\
\noindent Clear Water Bay, Kowloon, Hong Kong, China) \\
\smallskip

\noindent \textsc{B. Chen} \\
\medskip
\noindent E-mail: \url{mabfchen@ust.hk}\\
\noindent Department of Mathematics, \\
\noindent Hong Kong University of Science and Technology, \\
\noindent Clear Water Bay, Kowloon, Hong Kong, China \\
\smallskip

\noindent \textsc{H. Yang} \\
\noindent Department of Mathematics, \\
\noindent Guizhou University, Guiyang, 550025, China \\

\end{document}